\documentclass[10pt]{mcom-l}
\usepackage{amssymb}
\usepackage{amsmath}
\usepackage{microtype}
\usepackage{latexsym,enumerate,graphicx}
\usepackage{array,subfig}
\usepackage{hyperref}
\newtheorem{theorem}{Theorem}[section]
\newtheorem{lemma}[theorem]{Lemma}
\newtheorem{proposition}[theorem]{Proposition}
\newtheorem{corollary}[theorem]{Corollary}
\theoremstyle{definition}

\theoremstyle{remark}
\newtheorem{remark}{Remark}[section]

\newcommand{\supp}{\operatorname{supp}}

\newcommand{\dt}{\Delta t}
\newcommand{\ds}{\Delta s}
\newcommand{\dx}{\Delta x}
\newcommand{\A}{\ensuremath{\mathcal{A}}}
\newcommand{\B}{\ensuremath{\mathcal{B}}}
\newcommand{\eps}{\ensuremath{\varepsilon}}
\newcommand{\R}{\ensuremath{\mathbb{R}}}
\newcommand{\N}{\ensuremath{\mathbb{N}}}
\newcommand{\Z}{\ensuremath{\mathbb{Z}}}
\newcommand{\CS}{\mathcal{S}}
\newcommand{\BS}{\ensuremath{\mathbb{S}}}
\newcommand{\del}{\partial}
\newcommand{\ol}{\overline}
\newcommand{\ra}{\rightarrow}
\newcommand{\alp}{\alpha}

\newcommand{\M}{\mathcal{M}}

\newcommand{\tr}{\mathrm{tr}}

\newcommand{\bk}{\bar{k}}

\newcommand{\bmon}{B}
\newcommand{\CC}{\mathcal{C}}
\newcommand{\Dt}{\dt}
\newcommand{\Dx}{\dx}

\newcommand{\tf}{\tilde{f}}
\newcommand{\ts}{\tilde{\sigma}}
\newcommand{\tU}{\tilde{U}}
\newcommand{\tx}{\tilde{x}}
\newcommand{\ty}{\tilde{y}}
\newcommand{\tY}{\tilde{Y}}
\newcommand{\tL}{\tilde{L}}
\newcommand{\tc}{\tilde{c}}
\newcommand{\tb}{\tilde{b}}
\newcommand{\I}{\mathcal{I}_{\dx}}

\newcommand{\m}[2]{\overline{#1}^{\theta,#2}}
\newcommand{\mm}[1]{\overline{#1}^{\theta}}
\newcommand{\me}[1]{\overline{#1}^{\eps,\theta}}
\def\sym{\BS^N}

\numberwithin{equation}{section}

\newlength{\figwidth}
\begin{document}

\title[SL Schemes]{Semi-Lagrangian schemes for linear and fully
  non-linear diffusion equations}

\author[Debrabant]{Kristian Debrabant}
\address[Kristian Debrabant]{University of Southern Denmark, Department of Mathematics and Computer Science, Campusvej 55, 5230 Odense M, Denmark}
\email{debrabant\@@imada.sdu.dk}

\author[Jakobsen]{Espen R.~Jakobsen}
\address[Espen R. Jakobsen]{
    Department of Mathematical Sciences,
    Norwegian University of Science and Technology,
    7491 Trondheim, Norway}
\email{erj\@@math.ntnu.no}
\thanks{Jakobsen was supported by the Research
Council of Norway through the project ``Integro-PDEs: Numerical
methods, Analysis, and Applications to Finance''.}

\keywords{Monotone approximation schemes, difference-interpolation methods,
  stability, convergence, error bound, degenerate parabolic equations,
  Hamilton-Jacobi-Bellman equations, viscosity solution.}

\subjclass[2010]{
65M12, 
65M15, 
65M06, 
35K10, 
35K55, 
35K65, 
49L25, 
49L20. 
}
\date{}

\begin{abstract}
  For linear and fully non-linear diffusion equations of
  Bellman-Isaacs type, we introduce a class of approximation
  schemes based on differencing and interpolation. As opposed to classical
  numerical methods, these schemes work for general diffusions
  with coefficient matrices that may be non-diagonal dominant and
  arbitrarily degenerate.
  In general such schemes have to have a wide stencil.  Besides
  providing a unifying framework for several known first order
  accurate schemes, our class of schemes includes new first and
  higher order versions.
  The methods are easy to implement and more efficient than
  some other known schemes. We prove consistency and stability of the
  methods, and for the monotone first order methods, we prove
  convergence in the general case and robust error estimates in the
  convex case. The methods are extensively tested.
\end{abstract}
\maketitle
\section{Introduction}
\label{Sec:intro}
In this paper we introduce and analyze a class of
approximation schemes for fully non-linear diffusion equations of
Bellman-Isaacs type,
\begin{align}
\label{E}
u_t-\inf_{\alpha\in \A}\sup_{\beta \in \mathcal B}\Big\{L^{\alpha,\beta} [u](t,x) +
c^{\alpha,\beta}(t,x) u +
f^{\alpha,\beta} (t,x)\Big\}&=0&&\text{in}\quad Q_T,\\
u(0,x)&=g(x)&&\text{in}\quad\R^N,
\label{IV}
\end{align}
where $Q_T:=(0,T]\times\R^N$, $\A$ and $\B$ are complete metric spaces, and
\begin{align*}
&L^{\alpha,\beta} [u](t,x)= \tr [a^{\alpha,\beta}(t,x) D^2u(t,x)] +
b^{\alpha,\beta}(t,x) D
u(t,x).
\end{align*}
The coefficients
$a^{{\alpha,\beta}}=\frac12\sigma^{\alpha,\beta}\sigma^{\alpha,\beta\,\top}$,
$b^{{\alpha,\beta}}$,
$c^{{\alpha,\beta}}$, $f^{{\alpha,\beta}}$ and the initial data $g$ take values
respectively in
$\sym$, the space of $N \times N$ symmetric matrices, $\R^N$, $\R$,
$\R$, and $\R$. We will only assume that $a^{\alpha,\beta}$ is
positive semi-definite, thus the equation is allowed to degenerate and hence not have smooth solutions in general. Under suitable assumptions
(see Section \ref{Sec:WP}), the initial value problem
\eqref{E}-\eqref{IV} has a unique, bounded, H\"older continuous,
viscosity solution $u$.
This function is the upper or
lower value of a stochastic differential game, or, if $\mathcal A$ or
$\mathcal B$ is a singleton, the value
function of a finite horizon, optimal stochastic control problem \cite{YZ}.

We introduce a family of schemes that we call Semi-Lagrangian (SL)
schemes. They are a type difference-interpolation schemes and arise as
time-discretizations of the following semi-discrete equation
\[u_t-\inf_{\alpha\in \A}\sup_{\beta \in \mathcal
  B}\Big\{L^{\alpha,\beta}_k [\I u](t,x) +
c^{\alpha,\beta}(t,x) u +
f^{\alpha,\beta} (t,x)\Big\}=0\quad\text{in}\quad X_{\dx}\times(0,T),\]
where $L^{\alpha,\beta}_k$ is a monotone difference approximation of
$L^{\alpha,\beta}$ and $\I$ is an interpolation operator on
the spatial grid $X_{\dx}$. For more details see Section
\ref{Sec:Scheme}. Typically these schemes are first order accurate wide-stencil
schemes, and if the matrix $a^{\alp,\beta}$ is bad enough, the stencil
has to keep increasing as the grid is refined to have
convergence. They include as special cases monotone schemes from
\cite{F,M,CF:Appr,CL:MCM}, new versions of these schemes,
and a new non-monotone spatially second order accurate compact
stencil scheme. There are three main advantages of these
schemes: (i) they are easy to understand and implement, (ii) they are
faster than some alternative methods, and (iii) they are consistent
and, if $\I$ is monotone, monotone for every positive  semi-definite diffusion matrix
$a^{\alp,\beta}=\frac12\sigma^{\alp,\beta}\sigma^{\alp,\beta\,\top}$. The last
point is important because monotone methods are known to converge to
the correct solution \cite{BS:Conv}, while non-monotone methods need
not converge \cite{Ob} or can even converge to false solutions \cite{PFV}.

Classical finite difference approximations (FDMs) of \eqref{E} are not
monotone (of positive type) unless the matrix $a^{\alpha,\beta}$
satisfies additional assumptions like e.\,g.\ being diagonally
dominant \cite{KD:Book}. More general assumptions are given in
e.\,g.\ \cite{BZ,DK:ConstCoeff} but at the cost
of increased stencil length. In fact, Dong and Krylov \cite{DK:ConstCoeff}
proved that {\em no fixed stencil FDM} can approximate equations with
a second derivative term involving a general positive semi-definite
matrix function $a^{\alpha,\beta}$. Note that this type of result has
been known for a
long time, see e.\,g.\ \cite{MW,CL:MCM}. Some very simple examples of
such ``bad'' matrices are given by
\[\begin{pmatrix}x^2_1&\frac12x_1x_2\\
\frac12x_1x_2 & x_2^2\end{pmatrix}\ \text{in}\ [0,1]^2,\quad
\begin{pmatrix}\alp^2&\alp\beta\\
\alp\beta & \beta^2\end{pmatrix}\ \text{for}\ \mathcal
A=\mathcal B=[0,1],\quad I-\frac{(Du)(Du)^\top}{|Du|^2},
\]
and these types of matrices appear in Finance, Stochastic Control
Theory, and Mean Curvature Motion. The third example leads to
quasi-linear equations and will not be considered here, we refer instead
to \cite{CL:MCM}.

To obtain convergent or monotone methods for problems involving
non-diagonally dominant matrices, we know of two strategies: (i) The
classical method of rotating the coordinate system locally to
obtain diagonally dominant matrices
$a^{\alpha,\beta}$, see e.\,g.\ Section 5.4 in \cite{KD:Book}, and (ii)
the use of wide stencil methods. The two solutions seem to be somewhat
related, but at least the defining ideas and implementation are
different. Both ways lead to methods that have reduced order
compared to standard schemes for  diagonally dominant problems, but
the first strategy seems much more difficult to implement.

We mention that the schemes of \cite{MZ} for stationary Bellman equations have
much in common with the schemes in this paper. However the two types of
schemes and their derivation and error analysis are different in general.
Other related wide-stencil schemes are the method of \cite{BZ}, which
is not of SL type, and various SL like schemes for
other types of equations, e.\,g.\ the Mean
Curvature Motion equation \cite{CL:MCM},  Monge-Amp\`{e}re
equations \cite{Ob2}, and non-local Bellman equations \cite{CJ}. The
terminology SL schemes is already used for schemes for
transport equations, conservation laws, and first order
Hamilton-Jacobi equations. In the Hamilton-Jacobi setting, these
schemes go back to the 1983 paper \cite{CD} of Capuzzo-Dolcetta.

The rest of this paper is organized as follows. In the next section we
explain the notation and state a well-posedness and regularity result
for \eqref{E}--\eqref{IV}. The SL schemes are motivated and defined in
an abstract setting in Section \ref{Sec:Scheme}, and in Section
\ref{Sec:Anal} we prove that they are consistent, $L^\infty$-stable, and, if $\I$ is monotone, also monotone and convergent. We provide several examples of SL
schemes in Section \ref{Sec:Ex}, including the linear interpolation
SL (LISL) scheme. This  is the
basic example of this paper, a first order accurate wide
stencil scheme that can be defined on unstructured grids.
Higher order interpolation is not monotone. But for essentially monotone
solutions, we use third order monotonicity preserving cubic Hermite
interpolation \cite{FC,EJL} to obtain new schemes called monotonicity preserving
cubic interpolation SL (MPCSL) schemes in Section
\ref{Sec:PCSL}. These compact stencil schemes are second
order accurate in space and first or second order accurate in time.

In Section \ref{Sec:Disc}  we discuss various issues concerning the SL schemes. We compare the LISL scheme to the scheme of Bonnans-Zidani
\cite{BZ} and find that the LISL scheme is easier to understand and
implement and is faster in general.
We also explain that the SL schemes
can be interpreted as collocation methods for derivative free
equations, and as dynamic programming equations of discrete stochastic
differential games or optimal control problems.

 In Sections \ref{Sec:ErrBnd}, \ref{Sec:ErrUnif}
and Appendix \ref{Sec:Pf}, we derive robust error estimates for
monotone schemes for convex
equations, i.\,e.\ when $\I$ is monotone and \B\ is a singleton in \eqref{E}.
They are obtained through the regularization method of
Krylov \cite{Kr:00} and apply to degenerate equations and non-smooth
solutions.  Finally, in Section
\ref{Sec:Num} our methods are extensively tested. In particular we
find a first indication that the LISL and MPCSL schemes yield much
faster methods to solve the finance problem of \cite{BBMZ}.

\section{Notation and well-posedness}
\label{Sec:WP}

In this section we introduce notation and assumptions,
and give a well-posedness and comparison result for the initial value
problem \eqref{E} -- \eqref{IV}.

We denote by $\leq$ the component by component ordering in $\R^M$ and
the ordering in the sense of positive semi-definite matrices in $\sym$.
The symbols $\wedge$ and $\vee$ denote the minimum respectively the maximum.
By $|\cdot|$ we mean the Euclidean vector norm
in any $\R^p$ type space, e.\,g.\ if  $X\in\R^{N\times P}$ (an $N \times
P$ matrix), then
$|X|^2=\sum_{i,j}|X_{ij}|^2=\tr(XX^\top)$ where $X^\top$ is the
transpose of $X$.
If $w$ is a bounded function from some set $Q' \subset \ol Q_\infty$
into $\R$, $\R^M$, or $\R^{N\times P}$, we set
\[ |w|_{0} = \sup_{(t,y)\in Q'}\,|w(t,y)|, \qquad
 [w]_{\delta} = \sup_{(t,x)\neq
  (s,y)}\,\frac{|w(t,x)-w(s,y)|}{(|x-y|+|t-s|^{1/2})^\delta},\]
and $|w|_{\delta}=  |w|_{0} + [w]_{\delta}$ for any $\delta\in(0,1]$.
Let  $C_b (Q')$ and $\CC^{0,\delta}(Q')$, $\delta\in(0,1],$ denote
  respectively the space of bounded continuous functions on $Q'$ and
  the subset of $C_b (Q')$ in which the norm  $|\cdot|_\delta$ is finite.
Typically $Q'=Q_T$ or $Q'=\R^N$, and we will always suppress the
domain $Q'$ when writing norms.

For simplicity, we will use the following assumptions on the data of
\eqref{E}--\eqref{IV}:
\smallskip
\newlength{\eqlength}
\setlength{\eqlength}{\linewidth}
\addtolength{\eqlength}{-0.1\linewidth}
\begin{gather}\tag{A1}\label{A1}
\parbox[t]{\eqlength}{For any $\alp \in \A$ and $\beta\in\mathcal B$,
$a^{{\alpha,\beta}}=\frac12\sigma^{\alpha,\beta}\sigma^{\alpha,\beta\,\top}$ for some $N\times P$
matrix $\sigma^{{\alpha,\beta}}$. Moreover, there is a constant $K$ independent
of ${\alpha,\beta}$ such that}\\\nonumber
|g|_1+|\sigma^{{\alpha,\beta}}|_1+|b^{{\alpha,\beta}}|_1 + |c^{{\alpha,\beta}}|_1+|f^{{\alpha,\beta}}|_1 \leq
K.
\end{gather}
These assumptions are standard and ensure comparison and well-posedness of
\eqref{E}--\eqref{IV} in the class of bounded $x$-Lipschitz functions.
\begin{proposition}
\label{WPi}
Assume that \eqref{A1} holds. Then there exist a unique solution $u$ of
\eqref{E}--\eqref{IV} and a constant $C$ only depending on $T$ and $K$
from \eqref{A1} such that
\[|u|_1\leq C.\]

Furthermore, if $u_1$ and $u_2$ are sub- and supersolutions of
\eqref{E} satisfying $u_1(0,\cdot)\leq u_2(0,\cdot)$,
then $u_1\leq u_2$.
\end{proposition}
The proof is standard. Assumption \eqref{A1} can be relaxed in many ways,
e.\,g.\ using weighted norms, H\"older or uniform continuity,
etc. But in doing so, solutions can become unbounded and less smooth, and
the analysis becomes harder and less transparent. Therefore we will
not consider such extensions in this paper.

By solutions in this paper we always mean viscosity solutions, see
e.\,g.\ \cite{CIL:UG,YZ}.

\section{Definition of SL schemes}
\label{Sec:Scheme}

In this section we propose a class of approximation  schemes
for \eqref{E}--\eqref{IV} which we call Semi-Lagrangian or SL
schemes. This class includes (parabolic versions of) the
``control schemes'' of  Menaldi \cite{M} and Camilli and Falcone
\cite{CF:Appr} and the monotone schemes of Crandall and Lions
\cite{CL:MCM}. It also includes SL schemes for first
order Bellman equations \cite{CD,F}, and some new versions as discussed in Section \ref{Sec:Ex}. For a
motivation for the name, see Remark \ref{rem_name}. For the
time discretization we use a new generalized mid-point rule
that includes explicit, implicit, and a second order Crank-Nicolson like
approximations. Since the equation  is non-smooth in general, the
usual  way of defining a Crank-Nicolson scheme \cite{crank47apm} only gives a
first order accurate scheme in time.

The schemes are defined on a possibly unstructured family of grids
$\{G_{\dt,\dx}\}$,
\[G=G_{\dt,\dx}=\{(t_n,x_i)\}_{n\in\N_0,i\in\N}=\{t_n\}_{n\in\N_0}\times X_{\dx},\]
for $\dt,\dx>0$. Here $0=t_0<t_1<\dots<t_n<t_{n+1}$ satisfy
\[\max_n\Delta t_n\leq\dt \qquad\text{where}\qquad \dt_n=t_{n}-t_{n-1} ,\]
and $X_{\dx}=\{x_i\}_{i\in\N}$ is the set of vertices or nodes for a
non-degenerate polyhedral subdivision
$\mathcal{T}^{\dx}=\{T^{\dx}_j\}_{j\in\N}$ of $\R^N$, i.\,e., the polyhedrons $T^{\dx}_j$ satisfy
\begin{gather*}\hbox{int}(T^{\dx}_j\cap T^{\dx}_i)\underset{i\neq j}{=}\emptyset, \quad
\underset{j\in\N}{\bigcup}\, T^{\dx}_j=\R^N,\\\rho \dx\leq
\sup_{j\in\N}\{ \hbox{diam}\,B_{T^{\dx}_j}\}\leq \sup_{j\in\N}\{
\hbox{diam}\,T^{\dx}_j\}\leq
 \dx\end{gather*}
for some
$\rho\in(0,1)$, where $\hbox{diam}$ is the diameter of the set
 and $B_{T^{\dx}_j}$ is the greatest ball contained in $T^{\dx}_j$.

To motivate the numerical schemes, we write
$\sigma=(\sigma_1,\sigma_2,...,\sigma_m,...,\sigma_P)$ where
$\sigma_m$ is the $m$-th column of $\sigma$ and observe
that for $k>0$ and smooth functions $\phi$,
\begin{align*}
\frac12\tr [\sigma\sigma^\top
  D^2\phi(x)]&=\frac12\sum_{m=1}^P\tr [\sigma_m\sigma_m^\top D^2\phi(x)]\\
&
  =\sum_{m=1}^P\frac12\frac{\phi(x+k\sigma_m)-2\phi(x)+\phi(x-k\sigma_m)}{k^2}+\mathcal{O}(k^2),\\
b D\phi(x) &= \frac{\phi(x+k^2b)-\phi(x)}{k^2}+\mathcal{O}(k^2)\\
&=\frac12\frac{\phi(x+k^2b)-2\phi(x)+\phi(x+k^2b)}{k^2}+\mathcal{O}(k^2) .
\end{align*}
These approximations are monotone (of positive type) and the errors
  are bounded by $\frac1{48}P|\sigma|_0^4|D^4\phi|_0k^2$ and
  $\frac12|b|_0^2|D^2\phi|_0k^2$ respectively. To relate these approximations
  to a grid $G$, we replace $\phi$ by its interpolant $\I\phi$ on that
  grid and obtain
\begin{align*}
\frac12\tr [\sigma\sigma^\top
  D^2\phi(x)] &
  \approx\sum_{m=1}^P\frac12\frac{(\I\phi)(x+k\sigma_m)-2(\I\phi)(x)+(\I\phi)(x-k\sigma_m)}{k^2},\\
b D\phi(x) &\approx \frac12\frac{(\I\phi)(x+k^2b)-2(\I\phi)(x)+(\I\phi)(x+k^2b)}{k^2}.
\end{align*}
If the interpolation is monotone (positive) then the full discretization is
still monotone and represents a typical example of the discretizations
we consider below.

To construct the general scheme, we generalize the above
construction. Consider general finite difference approximations of the
differential operator $L^{\alpha,\beta}[\phi]$ in \eqref{E} defined as
\begin{align}
\label{BoAppr}
&L_k^{\alpha,\beta}[\phi](t,x):=
\sum_{i=1}^M\frac{\phi(t,x+y^{{\alpha,\beta},+}_{k,i}(t,x))-2\phi(t,x)+\phi(t,x+ y^{{\alpha,\beta},-}_{k,i}(t,x))}{2k^2},
\end{align}
for $k>0$ and some $M\geq1$, where for all smooth functions $\phi$,
\begin{align}
\label{consistL}
|L^{\alpha,\beta}_k[\phi]-L^{\alpha,\beta}[\phi]|\leq C(|D\phi|_0+\dots+|D^4\phi|_0)k^2.
\end{align}
This approximation and interpolation yield a
semi-discrete approximation of \eqref{E},
\begin{align*}
 U_t-\inf_{\alp\in\mathcal A}\sup_{\beta \in \mathcal B}\Big\{L_k^{\alpha,\beta}[\I U](t,x)+c^{\alpha,\beta}(t,x)U+f^{\alpha,\beta}(t,x)
\Big\}=0\quad
\text{in}\quad (0,T)\times X_{\dx},
\end{align*}
and the final scheme can then be found after discretizing in time using a parameter
$\theta\in[0,1]$,
\begin{align}
\label{FD}
&\delta_{\dt_n}U^n_i=\inf_{\alpha\in\mathcal{A}}\sup_{\beta \in \mathcal B}\Big\{
L_k^{\alpha,\beta}[\m{\I U}{n}_\cdot]^{n-1+\theta}_i+c^{{\alpha,\beta},n-1+\theta}_i\m{U}{n}_i
+f^{{\alpha,\beta},n-1+\theta}_i\Big\}
\end{align}
in $G$, where $U^n_i=U(t_n,x_i)$,
  $f^{{\alpha,\beta},n-1+\theta}_i= f^{\alpha,\beta}(t_{n-1}+\theta\dt_n,x_i)$, $\dots$ for $(t_n,x_i)\in
  G$,
\begin{gather*}
\delta_{\dt}\phi(t,x)=\frac{\phi(t,x)-\phi(t-\dt,x)}{\Dt},\qquad\text{and}
\qquad
\m{\phi}{n}_\cdot=(1-\theta)\phi^{n-1}_\cdot+\theta\phi^n_\cdot,\\
\m{\I\phi}{n}_\cdot=(1-\theta)\I\phi^{n-1}_\cdot+\theta\I\phi^n_\cdot.
\end{gather*}
 As initial conditions we take
\begin{align}
&U_i^0=g(x_i)\quad\text{in}\quad X_{\dx}.\label{FD_BC}
\end{align}
\begin{remark}
For the choices $\theta=0,1$, and $1/2$ the time discretization corresponds to respectively explicit Euler, implicit Euler and midpoint rule. For $\theta=1/2$, the full scheme can be seen as generalized Crank-Nicolson type discretization.
\end{remark}

\section{Analysis of SL schemes}
\label{Sec:Anal}

In this section we prove that the SL scheme \eqref{FD} is consistent
and $L^\infty$-stable, and in the case when the interpolation (and
hence the scheme) is monotone, we
present existence, uniqueness, and convergence results for the
schemes. Error estimates are given in
Section \ref{Sec:ErrBnd} for the monotone convex case.

For the approximation $L_k^{\alpha,\beta}$ and interpolation $\I$ we assume
\begin{align}
&\tag{Y1}\label{eq:Y1}\begin{cases}
&\displaystyle\sum_{i=1}^M [y^{{\alpha,\beta},+}_{k,i}+ y^{{\alpha,\beta},-}_{k,i}]=
  2k^2b^{\alpha,\beta}+\mathcal{O}(k^4),\\
&\displaystyle\sum_{i=1}^M [
y^{{\alpha,\beta},+}_{k,i}y^{{\alpha,\beta},+\,\top}_{k,i}+
  y^{{\alpha,\beta},-}_{k,i}y^{{\alpha,\beta},-\,\top}_{k,i}
  ]
= 2k^2
  \sigma^{\alpha,\beta}\sigma^{{\alpha,\beta}\,\top}+\mathcal{O}(k^4),\hspace{1.2cm}
  \\
&\displaystyle\sum_{i=1}^M [
y^{{\alpha,\beta},+}_{k,i,j_1}y^{{\alpha,\beta},+}_{k,i,j_2}y^{{\alpha,\beta},+}_{k,i,j_3}+
  y^{{\alpha,\beta},-}_{k,i,j_1}y^{{\alpha,\beta},-}_{k,i,j_2}y^{{\alpha,\beta},-}_{k,i,j_3}
]=\mathcal{O}(k^4),\\
&\displaystyle\sum_{i=1}^M [
y^{{\alpha,\beta},+}_{k,i,j_1}y^{{\alpha,\beta},+}_{k,i,j_2}y^{{\alpha,\beta},+}_{k,i,j_3}y^{{\alpha,\beta},+}_{k,i,j_4}+
  y^{{\alpha,\beta},-}_{k,i,j_1}y^{{\alpha,\beta},-}_{k,i,j_2}y^{{\alpha,\beta},-}_{k,i,j_3}y^{{\alpha,\beta},-}_{k,i,j_4}
]=\mathcal{O}(k^4),\\
\end{cases}\\
& \nonumber\quad\text{for all $j_1,j_2,j_3,j_4=1,2,\dots,N$ indicating
  components of the $y$-vectors.}\\
\tag{I1}\label{I1}
&\parbox[t]{\eqlength}{There are $K\geq0, r\in\N$ such
that for all smooth functions $\phi$,}\\
&\nonumber\qquad |(\I\phi)-\phi|_0\leq
K|D^r\phi|_0\dx^r.\\
\tag{I2}\label{I2}
&\parbox[t]{\eqlength}{There is a set of non-negative functions
  $\{w_{j}(x)\}_j$ such that}\\
\nonumber&\qquad(\I\phi)(x)=\sum_j\phi(x_j) w_{j}(x),\\
\nonumber&\text{and for all $i,j\in \N$,}\\
\nonumber&\quad
w_{j}(x)\geq0,\qquad w_{i}(x_j)=\delta_{ij},\qquad\text{and}\qquad\sum_i w_{i}(x) \equiv1.\\
\tag{I2'}\label{I2'}
&\parbox[t]{\eqlength}{Assumption (I2) holds, but $w_j=w_{\phi,j}$
is allowed to depend on $\phi$.}
\end{align}

Under assumption \eqref{eq:Y1}, a Taylor expansion shows that
$L_k^{\alpha,\beta}$ is a second order consistent approximation
satisfying \eqref{consistL}. If we assume also \eqref{I1}, it then follows that $L^{\alpha,\beta}_k[\I\phi]$ is a consistent
approximation of $L^{\alpha,\beta}[\phi]$ if $\frac{\dx^r}{k^2}\to0$. Indeed
\begin{align*}
&|L^{\alpha,\beta}_k[\I\phi]-L^{\alpha,\beta}[\phi]|\leq |L^{\alpha,\beta}_k[\I\phi]-L^{\alpha,\beta}_k[\phi]|+
|L^{\alpha,\beta}_k[\phi]-L^{\alpha,\beta}[\phi]|,
\end{align*}
where $|L^{\alpha,\beta}_k[\phi]-L^{\alpha,\beta}[\phi]|$ is estimated in
  \eqref{consistL}, and by \eqref{I1},
\begin{align*}
|L^{\alpha,\beta}_k[\I\phi]-L^{\alpha,\beta}_k[\phi]|\leq
C|D^r\phi|_0\frac{\dx^r}{k^2}.
\end{align*}

\begin{remark}
  Assumption \eqref{eq:Y1} is similar to the local consistency
  conditions used in \cite{KD:Book}. The $O(k^4)$ terms insure that
  the method is second order accurate as $k\ra0$. Convergence will
  still be achieved if we relax $O(k^4)$ to $o(k^2)$ as
  $k\ra0$.
\end{remark}

\begin{remark}
 An interpolation satisfying \eqref{I2'} is said to be {\em positive}
 and preserves positivity of the  data. Such an interpolation
 does not use (exact) derivatives to reconstruct the function
 $\phi$, and it may be a non-monotone and non-linear
 operator, as in the case of monotonicity preserving cubic
 interpolation (see Section \ref{Sec:PCSL}).

When \eqref{I2} holds and $w_{\phi,i}=w_i$ is independent of $\phi$,
the interpolation is a linear operator and {\em monotone} in the sense that
$U\leq V$ implies that $\I U\leq\I V$. The $w_j$'s form a basis and
the relation $\sum_i w_{\phi,i}(x)\equiv 1$ follows readily from the
other assumptions in \eqref{I1} and \eqref{I2}. Examples are constant,
linear, or multi-linear interpolation (i.\,e.\ $r\leq 2$ in
\eqref{I1}) since higher order interpolation is not monotone.
\end{remark}

The scheme is said to be of class $\tilde{L}$ if it can be written as
\begin{align}
\label{def_mon}
\sup_\alp\inf_\beta\Big\{
\bmon_{U^n,j,j}^{{\alpha,\beta},n,n} U_{j}^{n}
-\sum_{i\neq j}\bmon_{U^n,j,i}^{{\alpha,\beta},n,n}
U^{n}_i-\sum_{i} \bmon_{U^{n-1},j,i}^{{\alpha,\beta},n,n-1}
U_{i}^{n-1} - F^{{\alpha,\beta},n}_{j}\Big\} =0
\end{align}
in $G$, where $\bmon_{U^m,i,j}^{{\alpha,\beta},n,m}\geq 0$ might depend on $U^m$.
  In the case of monotone interpolation, $\bmon_{U^m,i,j}^{{\alpha,\beta},n,m}$ are independent of $U^m$, and the $\tilde{L}$-property implies monotonicity of the approximation scheme in the sense of Barles-Souganidis \cite{BS:Conv}.

 In the following, we denote by $c^{\alpha,\beta,+}$ the
positive part of $c^{\alpha,\beta}$. We now show consistency and
stability of the scheme.
\begin{lemma}[All SL schemes]
\label{LemConsistMon}
Assume \eqref{I1}, \eqref{I2'}, and \eqref{eq:Y1} hold.

\noindent (a) The scheme \eqref{FD} is consistent with \eqref{E} with
truncation error bounded by
\begin{align*}
\frac{|1-2\theta|}2|\phi_{tt}|_0\Dt
+C\Bigg(&\dt^2\left(
|\phi_{tt}|_0
+
|\phi_{ttt}|_0+
|D\phi_{tt}|_0+|
D^2\phi_{tt}|_0
\right)
\\
&+ |D^r\phi|_0
 \frac{\dx^r}{k^2}+(|D\phi|_0+\dots+|D^4\phi|_0)k^2
 \Bigg).
\end{align*}

\noindent (b) The scheme \eqref{FD} is of class $\tilde{L}$
(see \eqref{def_mon} for the definition) if
the following CFL condition holds,
\begin{align}
\label{CFL}
 (1-\theta)\Dt\Big[\frac
  {M}{k^2}-c^{{\alpha,\beta},n-1+\theta}_{i}\Big]\leq
  1\ \ \text{and}\ \ \theta\dt\,
  c^{{\alpha,\beta},n-1+\theta}_{i}\leq1 \ \text{for all}\ {\alpha,\beta},n,i.
\end{align}

\noindent (c) If in addition \eqref{A1} and \eqref{CFL} hold and
$2\theta\dt\sup_{\alpha,\beta}|c^{\alpha,\beta,+}|_0\leq 1$, then any solution $U$ of \eqref{FD}--\eqref{FD_BC} is $L^\infty$-stable satisfying
\begin{align*}
|U^n|_0&\leq e^{2\sup_{\alpha,\beta}|c^{\alpha,\beta,+}|_0t_n}\Big[|g|_0+t_n\sup_{\alpha,\beta}|f^{\alpha,\beta}|_0\Big].
\end{align*}
\end{lemma}

\begin{remark}
 By parabolic regularity $D^2\sim \partial_t$, so
  $|D^2\phi_{tt}|_0\sim|\phi_{ttt}|_0$.  When $\theta=1/2$, the
  scheme \eqref{FD} is second order accurate in time.
\end{remark}

\begin{proof}
(a) The scheme \eqref{FD} is consistent with
\eqref{E} with a truncation error bound
\begin{align*}
  \frac{|1-2\theta|}2|\phi_{tt}|_0\Dt+\frac13|\phi_{ttt}|_0\Dt^2
  +\sup_{\alp,\beta,n}\big\{\big|L^{\alpha,\beta}[\m{\phi}{n}]
  -L^{\alpha,\beta}_k[\m{\I \phi}{n}]\big|_0\big\}
  \\
  +\sup_{\alp,\beta,n}\Big\{\big|L^{\alpha,\beta}[\phi^{n-1+\theta}-\m{\phi}{n}]\big|_0+\big|
  c^{\alpha,\beta,n-1+\theta}(\phi^{n-1+\theta}-\m{\phi}{n})\big|_0
  \Big\}
\end{align*}
for smooth $\phi$. Part (a) now follows since by \eqref{I1},
\eqref{consistL},  and simple computations,
\begin{align*}
|L^{\alpha,\beta}[\m{\phi}{n}]-L^{\alpha,\beta}_k[\m{\I \phi}{n}]|&\leq C|D^r\phi|_0
 \frac{\dx^r}{k^2}+ C(|D\phi|_0+\dots+|D^4\phi|_0)k^2,\\
|L^{\alpha,\beta}[\phi^{n-1+\theta}-\m{\phi}{n}]|&\leq\dt^2\theta(1-\theta)
\sup_{\alp,\beta}|L^{\alpha,\beta}[\phi_{tt}]|_0\\
&\leq
C\dt^2 (|D\phi_{tt}|_0+|
D^2\phi_{tt}|_0),\\
|c^{\alpha,\beta,n-1+\theta}(\phi^{n-1+\theta}-\m{\phi}{n})\big|
&\leq C\theta(1-\theta)\dt^2|\phi_{tt}|_0.
\end{align*}
\medskip

\noindent (b) Note that since $\sum_i w_{\phi,i}\equiv 1$,
 \begin{gather*}
 L^{{\alpha,\beta}}_k[\I\phi(t,\cdot)](t_{n-1+\theta},x_j)=\sum_{i\in\N}
 l_{\phi,j,i}^{{\alpha,\beta},n-1+\theta}\big[\phi(t,x_i)-\phi(t,x_j)\big],
\intertext{where}
l_{\phi,j,i}^{{\alpha,\beta},n-1+\theta}=\sum_{l=1}^M
\frac{w_{\phi,i}\big(x_{j}+y_{k,l}^{{\alpha,\beta},+}(t_{n-1+\theta},x_{j})\big)+w_{\phi,i}\big(x_{j}+
  y_{k,l}^{{\alpha,\beta},-}(t_{n-1+\theta},x_{j})\big)}{2k^2}\end{gather*}
with $\sum_il_{\phi,j,i}^{{\alpha,\beta},n-1+\theta}=\frac {M} {k^2}$.
The $l_{\phi,j,i}^{{\alpha,\beta},n-1+\theta}$'s are non-negative by
\eqref{I2}.
The coefficients in \eqref{def_mon} can now be written as
\begin{align*}
\begin{array}{l}\bmon_{U^n,j,j}^{{\alpha,\beta},n,n}=1+\theta\Dt_n\,(\frac
  {M}
  {k^2}-l_{U^n,j,j}^{{\alpha,\beta},n-1+\theta}-c^{\alp,\beta,n-1+\theta}_{j}),
  \\
\\
\bmon_{U^{n-1},j,j}^{\alp,\beta,n,n-1}=1-(1-\theta)\Dt_n(\frac {M}
{k^2}-l_{U^{n-1},j,j}^{{\alpha,\beta},n-1+\theta}-c^{\alp,\beta,n-1+\theta}_{j}),
\\
\\
\bmon_{U^{n},j,i}^{{\alpha,\beta},n,n}=
 \theta\Dt_n
  l_{U^{n},j,i}^{\alp,\beta,n-1+\theta}, \qquad\qquad
  \bmon_{U^{n-1},j,i}^{\alp,\beta,n,n-1}=(1-\theta)\Dt_n
  l_{U^{n-1},j,i}^{\alp,\beta,n-1+\theta},
\end{array}
\end{align*}
where  $j\neq i$. These coefficients are positive if \eqref{CFL}
holds.
\medskip

\noindent (c) Fix any $\eps>0$ and let $j$ be such that
$|U^n_j|\geq |U|_0-\eps$. Assume first that $U^{n}_j\geq0$. By
the definition and sign of the $B$-coefficients (see part (b)),
\begin{align*}
&\bmon_{U^n,j,j}^{{\alpha,\beta},n,n}U^n_j \geq \Big(1-\theta\Dt_n\sup_{\alpha,\beta}|c^{\alpha,\beta,+}|_0\Big)U^{n}_j+\theta\Dt_n\Big(\frac  {M}
  {k^2}-l_{U^n,j,j}^{{\alpha,\beta},n-1+\theta}\Big)(|U^n|_0-\eps),\\
&-\sum_{i\neq j}\bmon_{U^n,j,i}^{{\alpha,\beta},n,n}
U^{n}_i\geq -\theta\Dt_n\Big(\frac  {M}
  {k^2}-l_{U^n,j,j}^{{\alpha,\beta},n-1+\theta}\Big)|U^n|_0, \\
&-\sum_{i} \bmon_{U^{n-1},j,i}^{{\alpha,\beta},n,n-1}
U_{i}^{n-1}\geq -\Big(1+(1- \theta)\Dt_n\sup_{\alpha,\beta}|c^{\alpha,\beta,+}|_0\Big)|U^{n-1}|_0.
\end{align*}
By \eqref{def_mon} we then find that
\begin{align*}
|U^n_j|&=U^n_j\\
&\leq
\Big(\frac{1+(1-\theta)\Dt\sup_{\alpha,\beta}|c^{\alpha,\beta,+}|_0}
{1-\theta\Dt\sup_{\alpha,\beta}|c^{\alpha,\beta,+}|_0}\Big)
\Big[|U^{n-1}|_0+\Dt_n\sup_{\alpha,\beta}|f^{\alpha,\beta}|_0+\theta\Dt_n\frac M{k^2}\eps\Big]\\
&\leq e^{2\sup_{\alpha,\beta}|c^{\alpha,\beta,+}|_0t_n}\Big[|g|_0+t_n\sup_{\alpha,\beta}|f^{\alpha,\beta}|_0\Big]+O_{k,n}(\eps).
\end{align*}
If $U^{n}_j<0$, a similar argument shows that $\tilde U=-U$
 satisfies the same inequality. Since
$|U^n|_0\leq |U^n_j|+\eps$ and $\eps$ is arbitrary, the result then follows.
\end{proof}

If monotone interpolation is used, we also prove existence,
uniqueness, and convergence of the schemes.
\begin{theorem}[Monotone SL schemes]
\label{FDthm}
Assume \eqref{A1}, \eqref{I1}, \eqref{I2}, \eqref{eq:Y1}, and \eqref{CFL}.
\smallskip

\noindent (a) There exists a unique bounded solution $U$ of
\eqref{FD}--\eqref{FD_BC}.
\smallskip

\noindent (b) $U$ converges uniformly to the solution $u$ of
\eqref{E}--\eqref{IV} as $\dt,k,\frac{\dx^r}{k^2}\rightarrow 0$.
\end{theorem}

\begin{proof}
(a) {\em Existence} and {\em uniqueness} follow by
induction. Let $t=t_n$ and assume $U^{n-1}$ is a known  bounded
function. For $\eps>0$ we define the operator $T$ by
\[T U^n_j = U^n_j -\eps\cdot (\text{left hand side of
  \eqref{def_mon}})\qquad\text{for all}\qquad j\in\Z^M.\]
Note that the fixed point equation $U^n=TU^n$ is equivalent to
  equation \eqref{FD}.
By the definition and sign of the $\bmon$-coefficients (which do not
depend on $U^n$ in this proof!) we see that
\begin{align*}
&TU^n_j-T\tU^n_j\\
&\leq
\sup_{\alp,\beta}\Big\{\Big[1-\eps\Big(1+\Dt_n\theta(\tfrac{M}{k^2}-l_{j,j}^{{\alpha,\beta},n-1+\theta}-c_j^{\alp,\beta,n-1+\theta})\Big)\Big](U^n_j-\tU^n_j)\\
&\qquad\quad+\eps\Dt_n\theta\Big(\tfrac{M}{k^2}-l_{j,j}^{{\alpha,\beta},n-1+\theta}\Big)|U^n_\cdot-\tU^n_\cdot|_0\Big\}\\
&\leq \Big(1-\eps\Big[1-\dt_n
\theta\sup_{\alpha,\beta}|c^{{\alpha,\beta},+}|_0\Big]\Big)|U^n_\cdot-\tU^n_\cdot|_0
\end{align*}
for $\eps$ such that
$1-\eps(1+\Dt\theta(\tfrac{M}{k^2}-c_j^{\alp,\beta,n-1+\theta}))\geq0$
and $\eps(1-\dt_n
\theta\sup_{\alpha,\beta}|c^{{\alpha,\beta},+}|_0)<1$ for all $j,n,{\alpha,\beta}$.
Taking the supremum over all $j$ and interchanging the role of $U$
and $\tU$ proves that $T$ is a contraction on
the Banach space of bounded functions on $X_{\dx}$ under the
  $\sup$-norm.
Existence and uniqueness then follows from the fixed
  point theorem (for $U^n$) and for all of $U$ by induction since
  $U^0=g$ is bounded.

\medskip
\noindent (b) In view of the $L^\infty$-stability of the scheme (Lemma
\ref{LemConsistMon} (c)), {\em convergence} of $U$ to the solution $u$
of \eqref{E}--\eqref{IV} follows from the Barles-Souganidis result in
\cite{BS:Conv}.
\end{proof}

\section{Examples of SL schemes}
\label{Sec:Ex}

\subsection{Examples of approximations \texorpdfstring{$L^{\alpha,\beta}_k$}{}}
\label{exLalp}
We present several examples of approximations of the
term $L^{\alpha,\beta}[\phi]$ of the form $L^{\alpha,\beta}_k[\phi]$, including previous
approximations that have appeared in \cite{F,M,CF:Appr,CL:MCM} plus
more computational efficient variants.

\begin{enumerate}[1.]
\item The approximation of Falcone \cite{F} (see also \cite{CD}),
\[b^{{\alpha,\beta}}D\phi \approx \frac{\I\phi(x+hb^{\alpha,\beta})-\I\phi(x)}{h},\]
corresponds to our $L^{\alpha,\beta}_k$ if $k=\sqrt h$, $y_k^{{\alpha,\beta},\pm}=k^2b^{\alpha,\beta}$.

\item\label{Appr2} The approximation of Crandall-Lions \cite{CL:MCM},
\[\frac12\tr[\sigma^{{\alpha,\beta}}\sigma^{{\alpha,\beta}\,\top}D^2\phi]\approx \sum_{j=1}^P\frac{\I\phi(x+k\sigma_j^{\alpha,\beta})-2\I\phi(x)+\I\phi(x-k\sigma_j^{\alpha,\beta})}{2k^2},\]
corresponds to our $L^{\alpha,\beta}_k$ if $y_{k,j}^{{\alpha,\beta},\pm}=
\pm k\sigma_j^{\alpha,\beta}$ and $M=P$.

\item The corrected version of the approximation of Camilli-Falcone
\cite{CF:Appr} (see also \cite{M}),
\begin{align*}
&\frac12\tr[\sigma^{{\alpha,\beta}}\sigma^{{\alpha,\beta}\,\top}D^2\phi]+b^{{\alpha,\beta}}D\phi\\
& \approx
\sum_{j=1}^P\frac{\I\phi(x+\sqrt{h}\sigma_j^{\alpha,\beta}+\frac{h}{P}b^{\alpha,\beta})-2\I\phi(x)+\I\phi(x-\sqrt{h}\sigma_j^{\alpha,\beta}+\frac{h}{P}b^{\alpha,\beta})}{2h},
\end{align*}
corresponds to our $L^{\alpha,\beta}_k$ if $k=\sqrt h$, $y_{k,j}^{{\alpha,\beta},\pm}=
\pm k\sigma_j^{\alpha,\beta}+\frac{k^2}{P}b^{\alpha,\beta}$  and $M=P$.

\item\label{Appr4} The new approximation obtained by combining approximations 1 and 2,
\begin{align*}
&\frac12\tr[\sigma^{{\alpha,\beta}}\sigma^{{\alpha,\beta}\,\top}D^2\phi]+b^{{\alpha,\beta}}D\phi\\
&\approx \frac{\I\phi(x+k^2b^{\alpha,\beta})-\I\phi(x)}{k^2}+
\sum_{j=1}^P\frac{\I\phi(x+k\sigma_j^{\alpha,\beta})-2\I\phi(x)
  +\I\phi(x-k\sigma_j^{\alpha,\beta})}{2k^2},
\end{align*}
corresponds to our $L^{\alpha,\beta}_k$ if $y_{k,j}^{{\alpha,\beta},\pm}= \pm
k\sigma_j^{\alpha,\beta}$ for $j\leq P$, $y_{k,P+1}^{{\alpha,\beta},\pm}= k^2 b^{\alpha,\beta}$ and $M=P+1$.
\item\label{Appr5} The new, more efficient version of approximation 3,
\begin{multline*}
\frac12\tr[\sigma^{{\alpha,\beta}}\sigma^{{\alpha,\beta}\,\top}D^2\phi]+b^{{\alpha,\beta}}D\phi
\approx
\sum_{j=1}^{P-1}\frac{\I\phi(x+k\sigma_j^{\alpha,\beta})-2\I\phi(x)+\I\phi(x-k\sigma_j^{\alpha,\beta})}{2k^2}\\
+\frac{\I\phi(x+k\sigma_P^{\alpha,\beta}+k^2b^{\alpha,\beta})-2\I\phi(x)+\I\phi(x-k\sigma_P^{\alpha,\beta}+k^2b^{\alpha,\beta})}{2k^2},
\end{multline*}
corresponds to our $L^{\alpha,\beta}_k$ if $y_{k,j}^{{\alpha,\beta},\pm}= \pm
k\sigma_j^{\alpha,\beta}$ for $j<P$, $y_{k,P}^{{\alpha,\beta},\pm}= \pm k\sigma_P^{\alpha,\beta}+k^2 b^{\alpha,\beta}$ and $M=P$.
\end{enumerate}

Approximation 5 is always more efficient than 3 in the
sense that it requires fewer arithmetic operations. The most
efficient of approximations 3, 4, and 5, is 4 when $\sigma^{\alpha,\beta}$ does not depend on ${\alpha,\beta}$ but $b^{\alpha,\beta}$ does, and 5 in the other cases.

\subsection{Linear interpolation SL scheme (LISL)}
To keep the scheme \eqref{FD} monotone, linear or
multi-linear interpolation is the most accurate interpolation
one can use in general. In this typical case we call the full scheme
\eqref{FD}--\eqref{FD_BC} the LISL scheme, and we will now summarize
the results of Section \ref{Sec:Anal} for this special case.
\begin{corollary}\label{cor:LISL}
Assume that \eqref{A1} and \eqref{eq:Y1} hold.
\smallskip

\noindent (a) The LISL scheme is monotone if the CFL
conditions \eqref{CFL} hold.
\smallskip

\noindent (b) The truncation error of the LISL scheme is
$O(|1-2\theta|\dt +\dt^2+ k^2 +\frac{\dx^2}{k^2} )$, so
it is first order accurate when
$k=O(\Dx^{1/2})$ and $\Dt=O(\dx)$, or if $\theta=\frac12$, $\Dt=O(\dx^{1/2})$.
\smallskip

\noindent (c)  If
$2\theta\dt\sup_{\alpha,\beta}|c^{\alpha,\beta,+}|_0\leq1$ and
\eqref{CFL} hold, then there exists a unique bounded and
$L^\infty$-stable solution $U$ of the LISL scheme converging
uniformly to the solution $u$ of \eqref{E}--\eqref{IV} as
$\dt,k,\frac{\dx}k\rightarrow 0$.
\end{corollary}

From this result it follows that the scheme is at most {\em first
  order accurate}, has {\em wide and increasing stencil} and a {\em
  good CFL condition}.  From the truncation error and the definition
of $L_k^{\alpha,\beta}$ the stencil is wide since the scheme is
consistent only if $\dx/k\ra0$ as $\dx\ra0$ and has stencil length
proportional to
  \[l:=\frac{\underset{t,x,{\alpha,\beta},i}{\max} |y^{{\alpha,\beta},-}_{k,i}|\vee
    |y^{{\alpha,\beta},+}_{k,i}|}{\dx}\sim \frac
  k{\dx}\ra\infty\quad\text{as}\quad\dx\ra0.\] Here we have used that
  if \eqref{A1} holds and $\sigma\not\equiv 0$, then typically
  $y^{{\alpha,\beta},\pm}_{k,i}\sim k$.  Note that if
  $k=\dx^{1/2}$, then $l\sim \dx^{-1/2}.$ Finally, in the case
  $\theta\neq1$ the CFL condition for \eqref{FD} is $\dt\leq
  Ck^2\sim\dx$ when $k=O(\dx^{1/2})$, and it is much less restrictive
  than the usual parabolic CFL condition, $\dt=O(\dx^2)$.

\subsection{A high order SL scheme for monotone solutions}
\label{Sec:PCSL}
In this section we introduce spatially second order accurate SL schemes
\eqref{FD}--\eqref{FD_BC} for non-degenerate tensor product grids. These
schemes are based on monotonicity preserving cubic (MPC) Hermite
interpolation \cite{FC,EJL} and will be denoted MPCSL schemes in
short. They are consistent for monotone (in
coordinate directions) solutions
of the scheme, but they are not monotone.

The MPC interpolation is obtained by a careful modification of
cubic Hermite interpolation \cite{EJL}, and for a function of
one variable on the interval $[x_i,x_{i+1}]$ it takes the form
\begin{gather}
\label{MPCI}
(\I\phi)(x)=\phi_i+(\phi_{i+1}-\phi_i)P_{i}(x)\\
\intertext{where}
P_{i}(x)=\alpha_i\frac{x-x_i}{\Delta
      x}+(3-\beta_i-2\alpha_i)\left(\frac{x-x_i}{\Delta
        x}\right)^2-(2-\alpha_i-\beta_i)\left(\frac{x-x_i}{\Delta
        x}\right)^3,
\end{gather}
where $\alp_i,\beta_i$ are bounded coefficients depending on
$\phi_{i-2},\phi_{i-2},\dots,\phi_{i+3}$. The algorithm is
described in Appendix \ref{Sec:MPCI}.
Multidimensional interpolation operators are obtained as
tensor products of one-dimensional interpolation operators, i.\,e.\ by
interpolating dimension by dimension.

\begin{remark}
\label{rem:cubint}
Rewriting $\I\phi$, we find that $(\I\phi)(x)=\sum_{i}\phi_iw_{\phi,i}(x)$ for
\[w_{\phi,i}(x)=(1-P_i(x))1_{[x_i,x_{i+1})}(x)+P_{i-1}(x)1_{[x_{i-1},x_{i})}(x)\]
and $1_{[x_i,x_{i+1}]}(x)$ is the indicator function that is $1$ in
$[x_i,x_{i+1})$ and $0$ otherwise. It is immediate that
$\sum_iw_{\phi,i}(x)\equiv 1$, and $w_{\phi,i}\geq0$ since
$P_i(x_i)=0$, $P_i(x_{i+1})=1$, and $P_i$ is monotone in between.
\end{remark}

\begin{lemma}
\label{propMC}
The above monotonicity preserving cubic interpolation satisfies
\eqref{I2'}. If the interpolated function is strictly monotone between
grid points, then \eqref{I1} holds with $r=4$ and the method is fourth
order accurate.
\end{lemma}
\begin{proof}
  Assumption \eqref{I2} holds by construction, see remark
  \ref{rem:cubint}. The
  error estimate follows from \cite{EJL}, since the above algorithm
  coincides with the two sweep algorithm given there when $n=1$
  interval is considered. In \cite{EJL} it is proved that this
  algorithm gives third order accurate approximations to the exact
  derivatives and hence the cubic Hermite polynomial constructed using
  this approximation is fourth order accurate.
\end{proof}

By Lemma \ref{propMC} and the results in Section \ref{Sec:Anal} we
have the following result:

\begin{corollary}
\label{propPCSL}
Assume \eqref{A1}, \eqref{eq:Y1} hold, and that for
all $\dx\in(0,1)$, solutions $U$ of the MPCSL scheme are such that
$\I U$ is strictly $x$-monotone between points in the $x$-grid $X_{\dx}$.
\smallskip

\noindent(a) The truncation error of the MPCSL scheme is
\[O\Big(|1-2\theta|\dt +\dt^2+k^2+\frac{\dx^4}{k^2}\Big),\]
and hence the scheme is second order accurate in space when $k=O(\dx)$
and first or second order accurate in time when $\theta\neq\frac12$
or $\theta=\frac12$ respectively.
\smallskip

\noindent(b) If $2\theta\dt\sup_{\alp,\beta}|c^{\alp,\beta}|_0\leq1$,
then the solution $U$ is $L^\infty$-stable.
\end{corollary}

\section{Discussion}
\label{Sec:Disc}
\subsection{Comparison with the scheme of Bonnans-Zidani (BZ)}
In \cite{BZ} (see also \cite{BOZ,BBMZ}) Bonnans and Zidani
suggest an alternative approach to discretize degenerate
diffusion equations. Their idea is to approximate the diffusion matrix
$a^{\alp,\beta}$ by a nicer matrix $a^{\alp,\beta}_k$ which admits monotone finite
difference approximations. For every $k\in\N$ they find a stencil
\[\CS_k\subset\{\xi=(\xi_1,\dots,\xi_N)\in\Z^N:~0<\max_{i=1}^N|\xi_i|\leq k,~i=1,\dots,N\}\] and positive numbers
$a^{\alpha,\beta}_{k,\xi}$ such that
\begin{align*}
a^{\alpha,\beta}\approx a^{\alp,\beta}_k:=\sum_{\xi\in\CS_k} a^{\alpha,\beta}_{k,\xi}
\xi\xi^\top.
\end{align*}
This leads to
a diffusion term that is a linear combination of directional
derivatives which are again approximated by central difference
approximations,
\begin{align*}
\tr[a^{\alp,\beta} D^2\phi] \approx \tr[a^{\alp,\beta}_k D^2\phi] =
\sum_{\xi\in\CS_k} a^{\alpha,\beta}_{k,\xi}D^2_\xi\phi\approx
\sum_{\xi\in\CS_k} a^{\alpha,\beta}_{k,\xi}\Delta_\xi\phi,
\end{align*}
where $D^2_\xi=\tr[\xi\xi^TD^2]=(\xi
\cdot D)^2$ and
\begin{align*}
&\Delta_{\xi} w(x)= \frac{1}{|\xi|^2\dx^2}\{w(x+\xi \dx) -2w(x)
  + w(x- \xi \dx)\}.
\end{align*}
This approximation is monotone by construction and respects the grid.
In two space dimensions, $a^{\alp,\beta}_k$ can be chosen such that
$|a^{\alp,\beta}-a^{\alp,\beta}_k|=O(k^{-2})$ (cf.\  \cite{BOZ}), and
then it is easy to see that the  truncation error is
\[O(k^{-2}+k^2\dx^2).\]

When $b^{\alp,\beta}\equiv0$, the BZ scheme can be obtained from \eqref{FD} by
replacing our $L^{\alp,\beta}_k$ by the above Bonnans-Zidani diffusion
approximation. This scheme shares many properties with the LISL scheme,
it is at most {\em first order accurate} (take $k\sim\dx^{-1/2}$), it has
a  similar {\em wide and  increasing  stencil}, and it has a similar {\em good
  CFL condition} $\dt\leq Ck^2\dx^2$ ($\sim \dx$ when
$k\sim\dx^{-1/2}$). To understand why the stencil is wide, simply note
that $k$ by definition is the stencil length and that the scheme is
consistent only if $k\ra\infty$ and $k\dx\ra0$.
The typical stencil length is $k\sim\dx^{-1/2}$, just as
it was for the LISL scheme.

The main drawback of this method is that it is costly since we must
compute the matrix $a^{\alp,\beta}_k$ for every $x,t,{\alp,\beta}$ in
the grid. In the fast two dimensional implementation in \cite{BOZ},
the number of operations for computing the coefficients for a fixed
$x,t,{\alp,\beta}$ is $\mathcal{O}(k)$ and thus goes to infinity as
$k\ra\infty$ in bad cases. The LISL scheme is easier to
understand and implement and is faster in the sense that the
computational cost for approximating the diffusion matrix for
fixed $x,t,{\alp,\beta}$ is independent of the stencil size. Later
we will see some numerical indication that the LISL scheme could be
faster than the BZ scheme in some test problems.

The MPCSL scheme in the typical case when $k=\dx$, is a {\em
  second order accurate in space} and {\em compact stencil} scheme having the usual
(not so good) {\em CFL conditions} for parabolic problems $\dt\sim
k^2=\dx^2$. When it can be used, it is far more efficient than the
other two schemes, see Section \ref{Sec:Num}. However there is no
proof that the method will converge to the correct solution, and it is
formally convergent only when the exact solutions are essentially
monotone, meaning monotone at least between grid points. Both the BZ
and LISL schemes ``always'' converge.
\subsection{Boundary conditions}
When solving PDEs on bounded domains, the SL (and BZ) schemes may
exceed the domain if they are not modified near the
boundary. The reason is of course the wide stencil. This may or may
not be a problem depending on the equation and the type of boundary
condition: (i) For Dirichlet conditions the scheme needs to be
modified near the boundary or boundary conditions must be
extrapolated. This may result in a loss of accuracy or monotonicity
near the boundary. (ii) Homogeneous Neumann conditions can be
implemented exactly by extending in the normal direction the values of
the solution on the boundary to the exterior.  (iii) If the boundary
has no regular points, no boundary conditions can be imposed. In this
case the SL schemes will not leave the domain if the normal
diffusion tends to zero fast enough when the boundary is
approached. Typical examples are equations of Black-Scholes type.

\subsection{Interpretation as a collocation method}
In the case the functions $w_{\phi,i}$ in \eqref{I2} do not
depend on $\phi$ (and form a basis), the scheme
\eqref{FD}--\eqref{FD_BC} can then be interpreted as a
collocation method for a derivative free equation, this is
essentially the approach of Falcone et al. \cite{F,CF:Appr}. The idea
is that if
\[W^{\dx}(Q_T)=\big\{u: u\text{ is a function on $Q_T$ satisfying } u\equiv
\I u \ \text{in}\ Q_T\big\}\]
denotes the interpolant space associated
to the interpolation $\I$, equation \eqref{FD} can be stated
in the following equivalent way: Find $U\in W^{\dx}(Q_T)$ solving
\begin{align}
\label{CE}
&\delta_{\dt_n}U^{n}_i=\inf_{\alpha\in\mathcal{A}}\sup_{\beta \in \mathcal B}\Big\{
L_k^{\alp,\beta}[\m{U}{n}]^{n-1+\theta}_i+c^{{\alpha,\beta},n-1+\theta}_i
\m{U}{n}_i
+f^{{\alpha,\beta},n-1+\theta}_i\Big\}\ \text{in}\ G.
\end{align}
In general $W^{\dx}$ can be any space of approximations which is
interpolating on the grid $X_{\dx}$, e.\,g.\  a space of splines, but we do not
consider this generality here.

\subsection{Stochastic game/control interpretation}
The scheme \eqref{FD}--\eqref{FD_BC} can be interpreted as
the dynamical programming equation of a discrete stochastic
differential game. We will explain this in the less
technical case when $\mathcal B$ is a singleton and the game
simplifies to an optimal stochastic control problem.

Assume that \eqref{A1} holds, and for simplicity, that
$c^{\alp}(t,x)\equiv0$ and the other coefficients are independent of $t$.
Then it is well-known (cf.\ \cite{YZ}) that the (viscosity) solution $u$ of
\eqref{E}--\eqref{IV} is the value function of the stochastic
control problem:
\begin{gather}
u(T-t,x)=\min_{\alp(\cdot)\in A}E\Big[\int_t^{T}
f^{\alp(s)}(X_s)\,ds+
g(X_T)\Big],\label{v_func}\\
\intertext{where $A$ is a set of admissible $\mathcal A$-valued
  controls and the diffusion process $X_s=X_s^{t,x,\alp(\cdot)}$
  satisfies the SDE}
 X_t=x \qquad\text{and}\qquad dX_s=\sigma^{\alp(s)}(X_s)\,dW_s+
 b^{\alp(s)}\,ds\quad \text{for}\quad s>t.\label{SDE}
 \end{gather}
This follows from dynamical programming (DP), and \eqref{E} is
called the DP equation for the control problem \eqref{v_func}--\eqref{SDE}.
Similarly, the schemes \eqref{FD}--\eqref{FD_BC} are DP equations (at
least in the explicit case)
of suitably chosen discrete time and space control problems approximating
\eqref{v_func}--\eqref{SDE}. We refer to \cite{KD:Book} for
more details.

We take the slightly
different approach explored in \cite{CD,F,M,CF:Appr} to show the
relation to control theory. The idea is to write the SL scheme in
collocation form \eqref{CE} and show that \eqref{CE} is the DP
equation of a discrete time continuous space optimal control problem.
We illustrate this approach by deriving an explicit scheme involving
$L_k^\alp$ as defined in part 4 Section \ref{exLalp}.
Let $\{t_0=0,t_1,\dots,t_M=T\}$ be discrete times and consider the
discrete time approximation of \eqref{v_func}--\eqref{SDE} given by
\begin{align}
&\tilde u(T-t_m,x)=\min_{\alp\in A_M}E\Big[\sum_{k=m}^{M-1}
f^{\alp_k}(\tilde X_k)\,\Delta t_{k+1}+
g(\tilde X_M)\Big],\label{v_func_approx}\\
&\tilde X_m=x,\quad
\tilde X_n=\tilde X_{n-1}+\sigma^{\alp_n}(\tilde X_{n-1})\,k_n\,\xi_n +
 b^{\alp_n}(\tilde X_{n-1})\,k_n^2\,\eta_n,\
 n>m,\label{SDEapprox}
\end{align}
where $k_n=\sqrt{(P+1)\Delta t_n}$, $A_M\subset A$ is an appropriate
subset of piecewise constant controls, and
$\xi_n=(\xi_{n,1},\dots,\xi_{n,P})^\top$ and $\eta_n$ are mutually
independent sequences of i.\,i.\,d.\  random variables satisfying
\begin{align*}P\Big((\xi_{n,1},\dots,\xi_{n,P},\eta_n)=\pm
e_j\Big)&= \frac1{2(P+1)}\quad\text{if } j\in\{1,\dots,P\},\\
P\Big((\xi_{n,1},\dots,\xi_{n,P},\eta_n)=e_{P+1}\Big)&=
\frac1{P+1},
\end{align*}
($e_j$ denotes the $j$-th unit vector) and all other values of $(\xi_{n,1},\dots,\xi_{n,P},\eta_n)$ have
probability zero. Here we have used a weak Euler approximation of the
SDE coupled with a quadrature approximation of the integral.
By DP
\[\tilde u(T-t_m,x)=\min_{\alp\in A_M}E\Big[\sum_{k=m}^{n-1}
f^{\alp_k}(\tilde X_k)\,\Delta t_{k+1}+
\tilde u(T-t_n,\tilde X_n)\Big] \quad\text{for all}\quad n>m,\]
and taking $n=m+1$, $s_{M-m}=T-t_{m}$, $\ds_m=s_{m}-s_{m-1}$, $\bk_m=k_{M-m}$, and evaluating the expectation using
\eqref{SDEapprox}, we see that
\begin{align*}
&\tilde u(s_{M-m},x)=\\
&\min_{\alp\in  \mathcal A}
\Big\{f^{\alp}(x)\ds_{M-m}+\frac{\bk_{M-m-1}^2}{P+1}L_{\bk_{M-m-1}}^{\alp}[\tilde
u](s_{M-m-1},x)+\tilde u(s_{M-m-1},x)\Big\},
\end{align*}
where $L_{k}^{\alp}$ is as in Section \ref{exLalp} part 4. If we subtract $\tilde u(s_{M-m-1},x)$ from
both sides and divide
by $\ds_{M-m}=
\frac{\bk_{M-m-1}^2}{P+1}$, we find \eqref{CE} with $\theta=0$.

In \cite{CF:Appr}, a similar argument is given in the stationary case
for schemes involving the $L^\alp_k$ of part 3 Section \ref{exLalp}. In
fact it is possible to identify all $L_k^\alp$'s appearing in Section
\ref{exLalp} with DP equations of suitably chosen discrete time
continuous space control problems. However assumption \eqref{eq:Y1}
is not strong enough for this approach to work for the general $L^\alp_k$
defined in Sections \ref{Sec:Scheme} and \ref{Sec:Anal}.

\begin{remark}
A DP approach naturally leads to explicit methods for time dependent
PDEs. But implicit methods can be derived from a
trick: Discretize the PDE in time by backward Euler to find a (sequence of)
stationary PDEs and use the DP approach on each stationary
PDE. This leads to an implicit iteration scheme since the DP equations
of stationary problems are always implicit.
\end{remark}

\begin{remark}
\label{rem_name}
By the definition of $L^\alp_k$ and \eqref{eq:Y1}, $x+y^{\alp,\pm}_{i,k}$ can be
seen as a short time approximation of \eqref{SDE}. Hence the scheme
\eqref{FD} tracks particle paths approximately. In view of the
discussion above we might say that the scheme follows particles in the mean
because of the expectation. For first order PDEs, schemes defined in
this way are called SL schemes by e.\,g.\ Falcone. Moreover, in this case our
schemes will coincide with the SL schemes of Falcone \cite{F} in the
explicit case. This explains why we choose to call these schemes SL
schemes also in the general case.
\end{remark}

\section{Error estimates in the monotone convex case}
\label{Sec:ErrBnd}
We derive error bounds when $\I$ is monotone and $\mathcal B$ is a
singleton and hence \eqref{E} is convex. It is not
known how to prove such results in the general case. In the following
we do not indicate the trivial $\beta$ dependence any more and we take
a uniform time-grid, $G=\dt\,\{0,1,\dots,N_T\}\times X_{\dx}$, for simplicity. Let
$Q_{\dt,T}:=\dt\,\{0,1,\dots,N_T\}\times\R^N$ and consider 
the intermediate equation
\begin{align}
\label{FD1}
&\delta_{\dt}V^n(x)
=\\
&\inf_{\alpha\in\mathcal{A}}\Big\{
  L_k^\alp[\m{V}{n}](t,x)+c^{{\alpha}}(t,x) \m{V}{n}(x)
+f^{{\alpha}}(t,x)\Big\}_{t=t_{n-1+\theta}}\quad
\text{in}\quad  \R^N\nonumber\\
\intertext{for $n=1,2,3,\dots$, with initial condition}
& V(0,x)= \ g(x)\quad\text{in}\quad \R^N.\label{FD1_BC}
\end{align}

\begin{lemma}
\label{ErrLem}
Assume that \eqref{I2} and the CFL condition \eqref{CFL} hold and that
$\sup_n |V^n|_1\leq C_V$. If $V$ solves \eqref{FD1}--\eqref{FD1_BC} and
$U$ solves \eqref{FD}--\eqref{FD_BC},
then
 \[|U-V| \leq C \frac{\dx}{k^2} \quad\text{in}\quad G.\]
\end{lemma}

\begin{proof}
Let $W=U-V$ and subtract the equation for $V$ from the one for
$U$ to find
\begin{align*}
W^n_i\leq &\  W^{n-1}_i+\dt\sup_{\alpha\in\mathcal{A}}\Big\{
L_k^{\alpha}[\m{\I W}{n}_\cdot]^{n-1+\theta}_i
+c^{{\alpha},n-1+\theta}_i\m{W}{n}_i\\
&+L_k^{\alpha}[\m{\I V}{n}-\m{V}{n}]_i^{n-1+\theta}\Big\}
\quad\text{in}\quad G.
\end{align*}
Let $C_c=\max_{\alp}|c^{\alp,+}|_0$. If $W_i^n\geq0$, we rearrange
using
\begin{align*}
|\I V^n-V^n|_0&\stackrel{\eqref{I2}}=\Big|\sum_{j}w_{j}(\cdot)\big(V_j^n-V^n(\cdot)\big)\Big|_0
\leq\Big|\sum_{j}w_{j}(\cdot)|V_j^n-V^n(\cdot)|\Big|_0
\\&
\leq\Dx|V^n|_1\Big|\sum_{j}w_{j}(\cdot)\Big|_0\stackrel{\eqref{I2}}=
\Dx|V^n|_1
\end{align*}
to see that
\begin{align*}
&\Big(1+\theta\dt\Big(\frac
{M}{k^2}-C_c\Big)\Big)W^n_i\\
&\leq  W^{n-1}_i
+ \dt\sup_{\alpha\in\mathcal{A}}\Big\{
    \theta
    \Big(L_k^{\alpha}[\I W^n_\cdot]^{n-1+\theta}_i +\frac {M}{k^2} W^n_i\Big)
\\&\quad
+(1-\theta)\Big(
L_k^{\alpha}[\I W^{n-1}_\cdot]^{n-1+\theta}_i+ c^{{\alpha},n-1+\theta}_iW^{n-1}_i
\Big)\Big\}\\
&\quad+ 2\dt\sup_n |V^n|_1\frac{\dx}{k^2}\quad\text{in}\quad G.
\end{align*}
By the CFL condition \eqref{CFL}, the coefficients of the above
inequality are all non-negative. Hence since $W^n\leq
|W^n|_0:=\sup_i|W^n_i|$, we may replace $W^n$ by $|W^n|_0$ on the
right hand side. Moreover,
since $\I|W^n|_0=|W^n|_0$ and $L^{\alp}_k[|W^n|_0]=0$, the upper bound
on the right hand side then reduces to
\[(1+\dt(1-\theta)C_c)|W^{n-1}|_0+\theta\dt\frac M{k^2}|W^n|_0+2C_V\dt\frac{\dx}{k^2}. \]
If $W_i^n<0$, then the same bound also holds for $-W_i^n$, and hence
\[(1+\dt\theta(\frac M{k^2}- C_c))|W^n|_0\leq (1+\dt(1-\theta)C_c)|W^{n-1}|_0+\theta\dt\frac M{k^2}|W^n|_0+2C_V\dt\frac{\dx}{k^2}\dt.\] Since $W^0\equiv 0$ in $X_{\dx}$, an
iteration then reveals that
\[|W^n|_0\leq
2C_V\dt\frac{\dx}{k^2}
\sum_{m=0}^n\Big(\frac{1+\dt(1-\theta)C_c}{1-\dt\theta C_c}\Big)^m\leq
t_n\frac{\dx}{k^2}4C_Ve^{C_ct_n}
\]
when $\dt$ is small enough, which implies the lemma.
\end{proof}

Next we estimate $|V-u|$, where $u$ solves \eqref{E}--\eqref{IV}, by the
regularization method of Krylov \cite{Kr:00}. To do that we need a continuity
and continuous dependence result for the scheme that relies on the following additional (covariance-type)
assumptions: Whenever two sets of data $\sigma, b$ and $\ts, \tb$ are
given, the corresponding approximations $L_k^\alp, y^{\alp,\pm}_{k,i}$
and $\tL_k^\alp, \ty^{\alp,\pm}_{k,i}$ in \eqref{BoAppr} satisfy
\begin{equation}\label{eq:Y3}
\tag{Y2}
\begin{cases}
  &\displaystyle\sum_{i=1}^M [y^{\alp,+}_{k,i}+
  y^{\alp,-}_{k,i}]-[\ty^{\alp,+}_{k,i}+ \ty^{\alp,-}_{k,i}]\leq
  2k^2(b^\alp-\tilde b^\alp),\\
  &\displaystyle\sum_{i=1}^M [
 y^{\alp,+}_{k,i}y^{\alp,+\,\top}_{k,i}+ y^{\alp,-}_{k,i}y^{\alp,-\,\top}_{k,i}
 ]+ [
  \ty^{\alp,+}_{k,i}\ty^{\alp,+\,\top}_{k,i}+\ty^{\alp,-}_{k,i}\ty^{\alp,-\,\top}_{k,i}
  ]\hspace{1cm}\\
  &\displaystyle\qquad -[
  y^{\alp,+}_{k,i}\ty^{\alp,+\,\top}_{k,i}+\ty^{\alp,+}_{k,i}y^{\alp,+\,\top}_{k,i}+y^{\alp,-}_{k,i}\ty^{\alp,-\,\top}_{k,i}+\ty^{\alp,-}_{k,i}y^{\alp,-\,\top}_{k,i}
] \\
  &\displaystyle \leq
  2k^2(\sigma^\alp-\ts^\alp)(\sigma^\alp-\ts^\alp)^\top +
  2k^4(b^\alp-\tb^\alp)(b^\alp-\tb^\alp)^\top,
\end{cases}
\end{equation}
when $\sigma,b,y_k^\pm$ are evaluated at $(t,x)$ and $\tilde\sigma,\tilde
b,\tilde y_k^\pm$ are evaluated at $(t,y)$ for all $t,x,y$.

In Section \ref{Sec:ErrUnif} we will prove the
following error estimate.
\begin{theorem}
\label{ErrBnd}
Assume that  $\mathcal B$ is a singleton, that \eqref{A1},
\eqref{eq:Y1}, \eqref{eq:Y3}, and the CFL conditions \eqref{CFL} hold,
and that $k\in(0,1)$ and $\dt\leq (2k_0\wedge 2k_1)^{-1}$. If $u$ and
$V$ are bounded solutions of \eqref{E}--\eqref{IV} and
\eqref{FD1}--\eqref{FD1_BC}, then
\[|V-u| \leq C(|1-2\theta|\dt^{1/4}+\dt^{1/3} +
k^{1/2})\quad\text{in}\quad Q_{\dt,T}.\]
\end{theorem}

It also follows from the regularity results in Section
\ref{Sec:ErrUnif} (see Proposition \ref{propV}) that $|V^n|_1\leq
2C_T$, so by Lemma \ref{ErrLem} and Theorem \ref{ErrBnd} we have
the following result.
\begin{corollary}[Error Bound]
Under \eqref{I1}, \eqref{I2}, and the assumptions of Theorem \ref{ErrBnd}, if $u$ solves
\eqref{E}--\eqref{IV} and $U$ solves \eqref{FD}--\eqref{FD_BC}, then
\[|u-U|\leq |u-V|+|V-U|\leq C(|1-2\theta|\dt^{1/4}+\dt^{1/3}+k^{1/2}+\frac{\dx}{k^2})\quad\text{in}\quad G.\]
\end{corollary}

This error bound applies to the LISL schemes,
and it also holds for unstructured grids. For more regular solutions
it is possible to obtain better error estimates, but general
and optimal results are not available.
The best estimate in our case is $O(\dx^{1/5})$ which is achieved when
$k=O(\dx^{2/5})$ and $\dt=O(k^2)$. Note that the CFL conditions
\eqref{CFL} already imply that $\dt=O(k^2)$ if
$\theta<1$.  Also note that the above bound does not show convergence
when $k$ is optimal for the LISL scheme ($k=O(\dx^{1/2})$).

\begin{remark}
These results are consistent with results for special
LISL type schemes for stationary Bellman equations. In fact if all
coefficients are independent of time and $c^\alp(x)<-c<0$, then by
combining the results of \cite{CF:Appr} and \cite{BJ:Rate}, exactly
the same error estimate is obtained for the solution of a particular
stationary LISL scheme and the unique {\em stationary} Lipschitz solution of
\eqref{E}.
\end{remark}

\section{Proof of Theorem \ref{ErrBnd}}
\label{Sec:ErrUnif}
We start by an existence and uniqueness result.

\begin{lemma}
\label{FD1thm}
Assume  that \eqref{A1}, \eqref{eq:Y1}, and
the CFL conditions \eqref{CFL} hold.
Then there exists a unique solution $U\in C_b(Q_{T,\dt})$ of
\eqref{FD1}--\eqref{FD1_BC}.
\end{lemma}
The  proof is similar to (but simpler than) the proof of Theorem
\ref{FDthm} with the modification that the fixed point is achieved in
the Banach space $C_b(\R^N)$ instead of the space of bounded
functions on $X_{\dx}$.

We now give a result comparing subsolutions of \eqref{FD1}
to supersolutions of
\begin{align}\label{FD2}
\begin{split}
\delta_{\dt}U^n(x)= &\inf_{\alpha\in\mathcal{A}}\Big\{
\tL_k^\alp[\m{U}{n}](t,x)+\tc^{{\alpha}}(t,x)\m{U}{n}
+\tf^{{\alpha}}(t,x)\Big\}_{t=t_{n-1+\theta}}
\text{in}\ \, \R^N, n\geq1,
\\U(0,x)=&\tilde g(x)\quad\text{in}\quad \R^N,
\end{split}
\end{align}
where $\tL_k^{\alpha}$ is the operator defined in \eqref{BoAppr},
\eqref{eq:Y1}, \eqref{eq:Y3} when $\sigma^\alp, b^\alp$ are replaced
by $\ts^\alp, \tb^\alp$.

\begin{theorem}
\label{LipCont} Assume
 that \eqref{A1}, \eqref{eq:Y3}, \eqref{CFL} hold for both \eqref{FD1}
and \eqref{FD2}. If $U\in
C(Q_{T,\dt})$ is a bounded above subsolution of
\eqref{FD1} and $\tU\in C(Q_{T,\dt})$ a bounded
below supersolution of \eqref{FD2}, then for all $k\in(0,1)$, $\dt\leq
(k_0\wedge k_1)^{-1}\wedge \frac{L_0}{2L}$ (see below), $x,y\in\R^N$, $n\in \{0,1,\dots,N_T\}$,
\begin{align*}
U(t_n,x)-\tU(t_n,y)\leq & \ R_{k_0}(t_n)|(U(0,\cdot)-\tU(0,\cdot))^+|_0 \\[0.1cm]
&+
  R_{k_0}(t_n)R_{k_1}(t_n)(L_0+t_nL)|x-y|
\\[0.1cm]
&+t_nR_{k_0}(t_n)\sup_{\alp\in\mathcal
  A}\big[|(f-\tf)^+|_0+(|U|_0\wedge|
    \tU|_0)|c-\tc|_0\big]\\[0.1cm]
&
+t_n^{1/2}2K_T\sup_{\alp\in\mathcal
  A}\big[|b-\tb|_0+|\sigma-\ts|_0\big]
\end{align*}
where
 $R_k(t)=1/(1-k\dt)^{t/\dt}$,
 $K_T \leq R_{k_0}(T)R_{k_1}(T)(L_0+TL)$,
\begin{align*}
&L_0=[g]_1\vee [\tilde g]_1+1, \quad L=([c^\alp]_1\vee[\tc^\alp]_1)
    (|U|_0\wedge |\tU|_0)+[f^\alp]_1\vee[\tf^\alp]_1,\\
& k_0=\sup_\alp|c^{\alp,+}|_0,\ \: \qquad
    k_1=8\sup_\alp\{[\sigma^\alp]_1^2+[b^\alp]_1^2+1\}.
\end{align*}
\end{theorem}

\begin{remark}
\label{remR}
The function $R_k(n\dt)=1/(1-k\dt)^n$ satisfies $\delta_{\dt}R_k(t_n)=k R_k(t_n)$,
$R_k(0)=1$, and $R_k(t_n)\leq e^{2kt_n}$ when $\dt\leq\frac1{2k}$.
\end{remark}

This is a key result in this paper, and the proof is given in Appendix
\ref{Sec:Pf}. In the stationary case,
results of this type have been obtained in \cite{BJ:Rate,CJ} for
simpler schemes. The result is a joint
uniqueness (take
$(\tilde\sigma,\tilde b,\tilde c, \tilde f, \tilde g)=(\sigma, b, c,
f, g)$), continuous dependence (take $x=y$), boundedness, and
$x$-Lipschitz continuity
result:
\begin{corollary}
\label{FD1apriori}
Under the assumptions of Theorem \ref{LipCont},
if $k\in(0,1)$
and $\dt\leq (2 k_0\wedge 2k_1)^{-1}$, then any bounded
solution $U\in C_b(Q_{T,\dt})$ of \eqref{FD1} satisfies
\smallskip

\noindent\,(i)\, $|U(t_n,\cdot)|_0\leq
e^{2k_0t_n}(|g|_0+t_n\sup_\alp|f^\alp|_0)$,
\smallskip

\noindent(ii) $|U(t_n,x)-U(t_n,y)|\leq e^{2(k_0+k_1)t_n}(L_0+t_nL)|x-y|,$
\smallskip

\noindent where the constants, which are defined in Theorem \ref{LipCont}, are
independent of $k,\dt,\dx$.
\end{corollary}

\begin{proof}
Part (i) follows from Theorem \ref{LipCont} (with $x=y$) and Remark \ref{remR}
since $\tU\equiv0$ satisfies \eqref{FD2} with $(\ts^\alp,\tb^\alp,\tc^\alp,\tf^\alp,\tilde g^\alp)=(\sigma^\alp, b^\alp, c^\alp, 0,0)$. Part (ii)
follows by taking $U=\tU$ and $x\neq y$.
\end{proof}

Now we extend the scheme \eqref{FD1} to the whole space $Q_T$.
One way to do this and to obtain continuous in time solutions is to
pose initial conditions on $[0,\dt)$ by interpolating between $g(x)$
and $U(\dt,x)$ where $U$ is the solution of \eqref{FD1}--\eqref{FD1_BC}.
\begin{align}
\label{FD11}
&\delta_{\dt}V(t,x)=\inf_{\alpha\in\mathcal{A}}\Big\{
  L_k^\alp[\mm{V}(t,\cdot)](t^\theta,x)+c^{\alpha}(t^\theta,x)\mm{V}(t,x)
+f^{\alpha}(t^\theta,x)\Big\}
\\ &\nonumber\hspace{8cm}\text{in}\  (\dt,T]\times\R^N,\\
&V(t,x)=\ \Big(1-\frac t{\dt}\Big)g(x)+\frac t{\dt}U(\dt,x)
\quad\text{in}\quad [0,\dt]\times\R^N.\label{FD11_BC}
\end{align}
where $\mm{V}(t,x)=(1-\theta)V(t-\dt,x) +\theta V(t,x)$ and
$t^\theta=t-(1-\theta)\dt$.
From the previous results for $U$ the existence, uniqueness, and
properties of $V$ easily follow.

\begin{proposition}
\label{propV}
Assume  that \eqref{A1}, \eqref{eq:Y1},
\eqref{eq:Y3}, and the CFL conditions \eqref{CFL} hold, and that
$k\in(0,1)$ and $\dt\leq (2k_0\wedge 2k_1)^{-1}$.
\smallskip

\noindent (a) There exists a unique solution $V\in
C_b(Q_T)$ of \eqref{FD11}--\eqref{FD11_BC}.
\smallskip

\noindent (b) There is a constant $C_T\geq0$ independent of $k,\dt,\dx$
such that
\begin{itemize}
\item [(i)\ ] $|V|_0\leq C_T$,\\[-0.3cm]
\item [(ii)\:] $|V(t,x)-V(t,y)|\leq C_T|x-y|\quad$ for all $\quad t\in [0,T],\
  x,y,\in\R^N,$\\[-0.3cm]
\item[(iii)] $|V(s_1,x)-V(s_2,x)|\leq C_T|s_1-s_2|^{1/2}\quad$ for all
  $\quad s_1,s_2\in [0,T],\  x,\in\R^N.$
\end{itemize}
\smallskip

\noindent (c) Let $V\in C_b(Q_T)$ and $\tilde V\in C_b(Q_T)$ be sub- and
supersolutions of \eqref{FD11}--\eqref{FD11_BC} corresponding to coefficients
$(\sigma^\alp,b^\alp,c^\alp,f^\alp,g)$ and $(\ts^\alp,\tb^\alp,\tc^\alp,\tf^\alp,\tilde g)$
respectively. Then there is a constant $C_T\geq0$ independent of $k,\dt,\dx$
such that for all $t\in[0,T]$,
\begin{align*}
|V(t,\cdot)-\tilde V(t,\cdot)|_0\leq C_T \Big(|g-\tilde
g|_0+t\sup_\alp[(|U|_0\wedge |\tU|_0)|c^\alp-\tc^\alp|_0+|f^\alp-\tf^\alp|_0]\\
+t^{1/2}\sup_\alp[|\sigma^\alp-\ts^\alp|_0+|b^\alp-\tb^\alp|_0]\Big).
\end{align*}
\end{proposition}

\begin{proof}
First note that the initial data on $[0,\dt]$ is uniformly bounded and
Lipschitz continuous in $x$ and $t$ by construction and Corollary
\ref{FD1apriori}.

(a) Existence of a bounded $x$-continuous solution follows from repeated
use of Lemma \ref{FD1thm} since we have initial conditions on
$[0,\dt]$. Continuity in time follows from Theorem \ref{LipCont} (with
$x=y$) since the data is $t$-continuous.

(b) Part (i) and (ii) follow from Corollary \ref{FD1apriori} since the
initial data is uniformly bounded and $x$-Lipschitz in $[0,\dt]$. To
prove part (iii) we assume $s_1<s_2$ and let $U(t,x)$ and $\tU(t,x)$ solve \eqref{FD11} with data
$
(\sigma^\alp(t+s_1,x),b^\alp(t+s_1,x),c^\alp(t+s_1,x),f^\alp(t+s_1,x),V(s_1,x))
\text{ and }(0,0,0,0,V(s_1,x))$
respectively. Note that for $t\in[0,T-s_1]$, $\tU(t,x)\equiv
V(s_1,x)$ and $U(t,x)\equiv V(t+s_1,x)$ where $V$ is the unique
solution of \eqref{FD11}--\eqref{FD11_BC}.
By part (c) we then get
\begin{multline*}
|V(t+s_1,\cdot)-V(s_1,\cdot)|_0=|U(t,\cdot)-\tU(t,\cdot)|_0\\
\leq C_T\Big(0 + t\sup_\alp[|f^\alp|_0+|V|_0|c^\alp|_0] +
  t^{1/2}\sup_\alp[|\sigma^\alp|_0+|b^\alp|_0]\Big)\quad\text{for}\quad t>0,
\end{multline*}
and hence part (iii) follows.

(c) Note that by construction of the initial data and Theorem
\ref{LipCont} with $x=y$, the result holds for $t\in[0,\dt]$, and then
the result holds for any $t>\dt$ by another application of Theorem
\ref{LipCont} with $x=y$.
\end{proof}

Using Krylov's method of shaking the coefficients \cite{Kr:00}, we will now find
smooth subsolutions of \eqref{FD11}.  First we introduce the
auxiliary equation
\begin{align}
    \label{FD11_shaken}
&\delta_{\dt}V^\eps(t,x)
=\inf_{\substack{0\leq s\leq \eps^2\\|e|\leq \eps\\
  \alpha\in\mathcal{A}}}
\Big\{
  L_k^\alp[\tau_{-e}\me{V}(t,\cdot)](r+s,x+e)\\
&\nonumber +c^{\alpha}(r+s,x+e) \me{V}(t,x)
+f^{\alpha}(r+s,x+e)\Big\}_{r=t^\theta-\dt-\eps^2}\
\text{in}\  (\dt,T]\times\R^N,\\
&V^\eps(t,x)= \Big(1-\frac t{\dt}\Big)g(x)+\frac t{\dt}V^\eps(\dt,x)
\quad\text{in}\quad [0,\dt]\times\R^N,\label{FD11_shaken_BC}
\end{align}
where $\tau_{e}\phi(t,x)=\phi(t,x+e)$ and $V^\eps(\dt,x)$ is
obtained by first solving \eqref{FD11_shaken} for discrete times $t_n=n\dt$.
For this equation to be well-defined for $t\in(\dt,T]$, the data and
$y^{\alp,\pm}_{k,i}$ must be defined for $t\in(-\dt-\eps^2,T+\eps^2]$. But this is ok
since  one can easily extend these functions to
$t\in[-r,T+r]$ for any $r>0$ in such a way that \eqref{A1},
\eqref{eq:Y1}, \eqref{eq:Y3} still hold. Also note that
\begin{multline}
\label{L-explain}
  L_k^\alp[\tau_{-e}\me{V}(t,\cdot)](r+s,x+e)
  =\ \frac1{2k^2}\sum_{i=1}^M\Big\{\me{V}(t,x+y^{\alp,+}_{k,i}(r+s,x+e))\\
  -2\me{V}(t,x)+\me{V}(t,x+y^{\alp,-}_{k,i}(r+s,x+e))
  \Big\},
\end{multline}
and hence \eqref{FD11_shaken} is an equation of the same type as
\eqref{FD11} (with different $\mathcal A$ and shifted
coefficients) satisfying \eqref{A1}, \eqref{eq:Y1}, \eqref{eq:Y3}
whenever \eqref{FD11} does.

By Proposition \ref{propV} there is a unique solution $V^\eps$ of
\eqref{FD11_shaken}--\eqref{FD11_shaken_BC} in
$[0,T+\dt+\eps^2]\times\R^N$. Let
$U^\eps(t,x):=V^\eps(t+\dt+\eps^2,x)$ and define by convolution,
\begin{align}
\label{Ueps}U_\eps(t,x)=\int_{\R^N}\int_0^\infty
U^\eps(t-s,x-e)\rho_\eps(s,e)\,ds\,de,
\end{align}
where $\eps>0$, $\rho_\eps(t,x)=\frac1{\eps^{N+2}}\rho(\frac
t{\eps^2},\frac x\eps)$, and
\[
\rho \in C^\infty(\R^{N+1}),\quad \rho\geq0,\quad \supp\,
\rho\subset[0,1]\times \{|x|\leq 1\}, \quad
\int_{\R^N}\rho(e) de =1.\]
Note that $U_\eps$ is well defined on the time interval $[-\dt,T]$. By
the next result it is the sought after smooth subsolution of \eqref{FD11}.

\begin{proposition}
\label{SubSol}
Under the assumptions of Proposition \ref{propV},
the function $U_\eps$ defined in \eqref{Ueps} satisfies
\medskip

\noindent (i) $U_\eps\in C^\infty((-\dt,T)\times\R^N)$, $|U_\eps|_1\leq C$,
$|D^m\del_t^n
U_\eps|_0\leq C\eps^{1-m-2n}$ for $n,m\in \N$.
\medskip

\noindent (ii) If $V$ is the solution of \eqref{FD11}--\eqref{FD11_BC}, then
 $|U_\eps-V|\leq C(\eps+\dt^{1/2})$ in $Q_T$.
\medskip

\noindent (iii) $U_\eps$ is a subsolution of \eqref{FD11} in $Q_T$.

\end{proposition}

\begin{proof}
The regularity estimates in (i) are immediate from properties of convolutions
and the regularity of $V^\eps$. The bound on $U_\eps-V$ (in $[0,T]$)
in (ii) follows from Proposition \ref{propV} (c) and \eqref{A1} which imply
\[|V^\eps-V|_0\leq C(\eps+\dt^{1/2}),\]
and regularity of $V^\eps$ along with properties of convolutions,
\[|U_\eps-V^\eps|_0\leq |U_\eps-U^\eps|_0+|V^\eps(\cdot+\dt+\eps^2,\cdot)-V^\eps|_0\leq
|V_\eps|_1(\eps+\dt^{1/2}).\]
To see that $U_\eps$ is a subsolution of \eqref{FD11}, first note
that from the definition of $U^\eps$ and \eqref{FD11_shaken} it follows that
\begin{multline*}
\delta_{\dt}U^\eps(t,x)\leq
  L_k^{\alp}[\tau_{-e}\me{U}(t,\cdot)](t^\theta+s,x+e)\\
\quad+c^{\alpha}(t^\theta+s,x+e)\me{U}(t,x)
+f^{\alpha}(t^\theta+s,x+e)
\end{multline*}
for all  $(t,x)\in[-\eps^2,T]\times\R^N$, $|e|,s^2\leq \eps$, and
$\alp\in\mathcal A$.
Now we change variables from $(t+s,x+e)$ to $(t,x)$ to find that
\begin{multline*}
\delta_{\dt}U^\eps(t-s,x-e)\leq
  L_k^{\alp}[\tau_{-e}\me{U}(t-s,\cdot)](t^\theta,x)\\
\quad+c^{\alpha}(t^\theta,x)\me{U}(t-s,x-e)
+f^{\alpha}(t^\theta,x)
\end{multline*}
for all  $(t,x)\in[0,T]\times\R^N$, $|e|,s^2\leq \eps$, and
$\alp\in\mathcal A$.
Then we multiply by $\rho_\eps(s,e)$ and integrate w.\,r.\,t.\ $(s,e)$. To
see what the result is, note that
\begin{multline*}
L_k^\alp[\tau_{-e}U^\eps(t-s,\cdot)](r,x)=\frac1{2k^2}\sum_{i=1}^M\Big\{
U^\eps(t-s,x+y^{\alp,+}_{k,i}(r,x)-e)\\
-2U^\eps(t-s,x-e)+U^\eps(t-s,x+y^{\alp,-}_{k,i}(r,x)-e)\Big\},
\end{multline*}
and hence
\begin{align*}
&\int\int
L_k^{\alp}[\tau_{-e}U^\eps(t-s,\cdot)](r,x)\rho_\eps(s,e)\,ds\,de
=  L^\alp_k[U_\eps(t,\cdot)](r,x).
\end{align*}
For the whole equation we then have
\begin{align*}
\delta_{\dt}U_\eps(t,x)&\leq
  L_k^{\alp}[\mm{U}_\eps(t,\cdot)](t^\theta,x)
  +c^{\alpha}(t^\theta,x)
  \mm{U}_\eps(t,x)
+f^{\alpha}(t^\theta,x)
\end{align*}
for all $(t,x)\in Q_T$ and $\alp\in\mathcal{A}$. Since this inequality
holds for all $\alp$, it follows that $U_\eps$ is a
subsolution of \eqref{FD11} in all of $Q_T$.
\end{proof}

We are now in a position to prove the error estimate
given in Theorem \ref{ErrBnd}.

\begin{proof}[Proof of Theorem \ref{ErrBnd}]
Let $U_\eps$ be defined in \eqref{Ueps}. By Proposition
\ref{SubSol} (i) and Lemma \ref{LemConsistMon} (a),
\begin{align*}
&\del_t U_\eps- \inf_{\alp\in\mathcal{A}}\Big\{L^\alpha
[\mm{U}_\eps(t,\cdot)](t^\theta,x)
+ c^\alpha(t^\theta,x) \mm{U}_\eps(t,x) +
f^\alpha (t^\theta,x)\Big\}\\
&\leq
\frac{|1-2\theta|}2|\del_t^2U_\eps|_0\Dt+
C\Big\{(|\del_t^2U_\eps|_0+|\del_t^3U_\eps|_0+|\del_t^2DU_\eps|_0+| \del_t^2D^2U_\eps|_0)\dt^2\\
&\quad+(|DU_\eps|_0+\dots+|D^4U_\eps|_0)k^2\Big\}\\
&\leq C\Big\{|1-2\theta|\eps^{-3}\dt+\eps^{-5}\dt^2+\eps^{-3}k^2\Big\}
\end{align*}
in $Q_T$. Moreover, by Proposition \ref{SubSol} (ii),
\[g(x)=U(0,x)\geq U_\eps(0,x)-C(\eps+\dt^{1/2}).\]
It follows that there is a constant $C\geq0$ such that
\[U_{\eps}-Ce^{\sup_\alp|c^\alp|_0t}\Big\{\eps+\dt^{1/2}
+t\Big(|1-2\theta|\eps^{-3}\dt+\eps^{-5}\dt^2
+\eps^{-3}k^2\Big)\Big\}\]
is a classical subsolution of \eqref{E}--\eqref{IV} with time shifted
coefficients. By continuous dependence and the comparison principle
\[
U_{\eps}-Ce^{\sup_\alp|c^\alp|_0t}\Big\{\eps+\dt^{1/2}
+t\Big(|1-2\theta|\eps^{-3}\dt+\eps^{-5}\dt^2
+\eps^{-3}k^2\Big)\Big\}\leq
u\ \ \text{in}\ \ Q_T,\]
and hence by Proposition \ref{SubSol} (ii),
\[U-u = (U-U_\eps)+(U_\eps-u)\leq
C\Big\{\eps+\dt^{1/2}+|1-2\theta|\eps^{-3}\dt+\eps^{-5}\dt^2+
\eps^{-3}k^2\Big\}.\]
We minimize w.\,r.\,t.\ $\eps$ and find that
\[u-U\leq
\begin{cases}
C (\dt^{1/4}+k^{1/2})&\text{ if }\theta\neq\frac12\\
C (\dt^{1/3}+k^{1/2})&\text{ if }\theta=\frac12
\end{cases}
\quad\text{in}\quad Q_T.\]

The lower bound on $u-U$ follows with symmetric -- but much easier -- arguments
where a smooth supersolution of the equation \eqref{E} is
constructed. Consistency and comparison for the scheme \eqref{FD11} is
then used to conclude. In view of Lemma \ref{LemConsistMon}, the
lower bound is a direct consequence of Theorem 3.1 (a) in \cite{BJ:Rate3}.
\end{proof}
\def\ignore#1\endignore{}
\newcolumntype{I}{@{}>{\ignore}c<{\endignore}}
\providecommand{\errtab}[3]{\begin{table}[ht]
\begin{tabular}{c|cIIcII|cIIcIIl}
&\multicolumn{6}{c|}{LISL}&\multicolumn{6}{c}{PCSL}\\ \cline{2-14}
$\dx$&error&&&rate&&&error&&&rate&&&
\\\hline
#1
\end{tabular}
\caption{#2}\label{#3}
\end{table}}
\providecommand{\errtimetab}[3]{\begin{table}
\begin{tabular}{c|cII|cII|c}
&\multicolumn{3}{c|}{error}&\multicolumn{3}{c|}{rate}&time in s\\\hline
$\dx$&$L^\infty$&$L^2$&$L^1$&$L^\infty$&$L^2$&$L^1$&\\\hline
#1
\end{tabular}
\caption{#2}\label{#3}
\end{table}}
\providecommand{\errtimesubtab}[3]{
\begin{tabular}[t]{c|cII|cII|r}
$\dx$&error&$L^2$&$L^1$&rate&$L^2$&$L^1$&time in s\\\hline
#1
\end{tabular}
}
\section{Numerical results}\label{Sec:Num}
In the following, we apply the LISL and MPCSL schemes to
linear and convex test problems in two space-dimensions, and hence have
no dependence of $\beta$.
For the LISL scheme, we choose $k=\sqrt{\dx}$ and a regular triangular
grid, whereas for the MPCSL scheme we choose $k=\dx$ and a regular
rectangular grid. If not stated otherwise, we use $\theta=0$ (explicit
methods), CFL condition $\Delta t=k^2$, and approximation
\ref{exLalp}.\ref{Appr5} for $L^{\alpha,\beta}$.
As error measure we will always use
the $L^\infty$-norm, and the error rates are calculated as $
r_i=\frac{\ln\|e_i\|-\ln\|e_{i-1}\|}{\ln\|\dx_i\|-\ln\|\dx_{i-1}\|}.$
All calculations are done in Matlab, on an INTEL(R) Core(TM)2 Duo P8700, 2.54Ghz Laptop.

\subsection{Linear problem with smooth solution}
Our first problem is taken from \cite{BOZ} and has exact solution
$u(t,x)=(2-t)\sin x_1\sin x_2$, its coefficients in \eqref{E} are
\begin{align*}
f^\alpha(t,x)=&\sin x_1\sin x_2[(1+2\beta^2)(2-t)-1]\\&-2(2-t)\cos x_1\cos x_2\sin(x_1+x_2)\cos(x_1+x_2),\\
c^\alpha(t,x)=&0,\qquad b^\alpha(t,x)=0,\qquad\sigma^\alpha(t,x)=\sqrt2\begin{pmatrix}
\sin(x_1+x_2)&\beta&0\\
\cos(x_1+x_2)&0&\beta
\end{pmatrix}.
\end{align*}
We consider $\beta^2=0.1$ and $\beta=0$. Note that in the second case,
the scheme considered in \cite{BOZ} is not consistent. Table \ref{Tab:Ex1aundb} gives the (spatial) errors and rates obtained
at $t=1$ applying the LISL and the MPCSL scheme, as well as the CPU time needed. As the solution is linear in $t$, one time step suffices.
\newcommand{\newtab}[3]{\begin{table}[ht]
\begin{tabular}{c|crr|crr}
$\dx$&\multicolumn{3}{c|}{LISL}&\multicolumn{2}{c}{MPCSL}\\ \cline{2-7}
&error&rate&time in s&error&rate&time in s\\\hline
#1
\end{tabular}
\caption{#2}\label{#3}
\end{table}
}
\newcommand{\newsubtab}[3]{
\subfloat[#2]{
\begin{tabular}{c|crr|crr}
$\dx$&\multicolumn{3}{c|}{LISL}&\multicolumn{2}{c}{MPCSL}\\ \cline{2-7}
&error&rate&time in s&error&rate&time in s\\\hline
#1
\end{tabular}
}}
\begin{table}[ht]
\newsubtab{
3.93e-2&3.79e-2&&0.07&1.03e-3&&0.11\\
1.96e-2&1.93e-2&0.97&0.30&2.57e-4&2.00&0.54\\
9.82e-3&9.45e-3&1.03&1.53&6.42e-5&2.00&3.75\\
4.91e-3&4.50e-3&1.07&6.25&1.61e-5&2.00&15.16\\
2.45e-3&2.43e-3&0.89&24.77&4.01e-6&2.00&62.02}
{$\beta^2=0.1$}{Tab:Ex1a}
\\\newsubtab{
3.93e-2&3.94e-2&&0.04&1.03e-3&&0.06\\
1.96e-2&1.98e-2&0.99&0.14&2.57e-4&2.00&0.24\\
9.82e-3&9.94e-3&0.99&0.68&6.43e-5&2.00&1.47\\
4.91e-3&4.70e-3&1.08&2.61&1.61e-5&2.00&5.94\\
2.45e-3&2.45e-3&0.94&10.64&4.02e-6&2.00&25.29}
{$\beta=0$}{Tab:Ex1b}
\caption{Results for the smooth linear problem at $t=1$, grid adapted to monotonicity\label{Tab:Ex1aundb}}
\end{table}

As expected for smooth solutions, in both cases we obtain order one
for the LISL scheme and order two for the MPCSL scheme, and the CPU
time needed is proportional to the number of grid points
$\frac1{\Delta x^2}$. Here, we
have chosen the grid points such that the solution is monotone in between.
If not, we would get order one for the LISL scheme but
no convergence for the MPCSL scheme (see Section \ref{Sec:PCSL}).

\subsection{Linear problem with non-smooth solution}
The second problem we test has a non-smooth
exact solution
\begin{equation*}
u(t,x)=(1+t)\sin\frac{x_2}2
\begin{cases}
\sin\frac{x_1}2&\text{ for }-\pi<x_1<0,\\
\sin\frac{x_1}4&\text{ for }0<x_1<\pi
\end{cases}
\end{equation*}
in $[-\pi,\pi]^2$ and coefficients in \eqref{E} given by
\begin{align*}
f^\alpha(t,x)=& \
\sin\frac{x_2}2
\begin{cases}
\sin\frac{x_1}2\left(1+\frac{1+t}{4}
(\sin^2x_1+\sin^2x_2)
\right)
&\text{ for }-\pi<x_1<0
\\
\sin\frac{x_1}4\left(
1+\frac{1+t}{16}
(
\sin^2x_1+4\sin^2x_2
)
\right)
&\text{ for }0<x_1<\pi
\end{cases}
\\&
-\sin x_1\sin x_2\cos\frac{x_2}2
\begin{cases}
\frac{1+t}2\cos\frac{x_1}2&\text{ for }-\pi<x_1<0
\\
\frac{1+t}4\cos\frac{x_1}4&\text{ for }0<x_1<\pi
\end{cases}
,\\
c^\alpha(t,x)=&0,\qquad
b^\alpha(t,x)=0,\qquad\sigma^\alpha(t,x)=\sqrt2\begin{pmatrix}
\sin x_1\\
\sin x_2
\end{pmatrix},
\end{align*}
and we pose Dirichlet boundary conditions. This is a monotone
non-smooth problem, and we obtain order one half applying the LISL
scheme and order one applying the MPCSL scheme, i.\,e.\ reduced rates,
see Table \ref{Tab:Ex2}. Again, one time step suffices, and the CPU
time needed is thus proportional to $\frac1{\Delta x^2}$.
\newtab{
7.76e-2&1.24e-2&&0.02&7.56e-3&&0.03\\
3.90e-2&8.75e-3&0.51&0.04&4.19e-3&0.86&0.06\\
1.96e-2&6.19e-3&0.50&0.14&2.20e-3&0.93&0.28\\
9.80e-3&4.38e-3&0.50&0.76&1.12e-3&0.97&1.79\\
4.90e-3&3.10e-3&0.50&2.99&5.69e-4&0.98&7.16\\
2.45e-3&2.19e-3&0.50&11.48&2.86e-4&0.99&28.66}
{Results for the non-smooth linear problem at $t=1$}{Tab:Ex2}

\subsection{Optimal control problems with smooth solutions}
\begin{enumerate}[(A)]
\item\label{OptProb1}
We test an example from \cite{BOZ} with exact solution
$u(t,x_1,x_2)=\left(\frac32-t\right)\sin x_1\sin x_2$. The
corresponding coefficients and control set in \eqref{E} are
\begin{align*}
f^\alpha=&\left(\frac12-t\right)\sin x_1\sin x_2
+\left(\frac32-t\right)\Bigg[
\sqrt{\cos^2x_1\sin^2x_2+\sin^2x_1\cos^2x_2}\\
&-2\sin(x_1+x_2)\cos(x_1+x_2)\cos x_1\cos x_2\Bigg],\\
c^\alpha=&\ 0,\quad
b^\alpha=\alpha,\quad\sigma^\alpha=\sqrt2\begin{pmatrix}
\sin(x_1+x_2)\\
\cos(x_1+x_2)
\end{pmatrix},\quad \A=\{\alpha\in\R^2:~\alpha_1^2+\alpha_2^2=1\}.
\end{align*}
 As $\sigma^\alpha$ does not depend on $\alpha$ but $b^\alpha$ does,
 we choose approximation \ref{exLalp}.\ref{Appr4} for
 $L^{\alpha,\beta}$ and thus need only about half of the number of
 interpolations we would need if we had chosen approximation
 \ref{exLalp}.\ref{Appr5}.\\

\item\label{OptProb2}
The next test problem has exact solution
$u(t,x_1,x_2)=\left(2-t\right)\sin x_1\sin x_2$ and coefficients and
control set given by
\begin{align*}
f^\alpha(t,x)=&\ (1-t)\sin x_1\sin x_2-2\alpha_1\alpha_2(2-t)\cos x_1\cos x_2,\\
c^\alpha(t,x)=&0,\quad
b^\alpha(t,x,\alpha)=0,\quad\sigma^\alpha=\sqrt2\begin{pmatrix}
\alpha_1\\
\alpha_2
\end{pmatrix},\quad\A=\{\alpha\in\R^2:~\alpha_1^2+\alpha_2^2=1\}.
\end{align*}
\end{enumerate}
In both examples, due to the solution being linear in $t$, one time
step suffices. The results at $t=0.5$ are given in Table
\ref{Tab:Ex3aundb}, where again the grid is adapted to
monotonicity. As expected for smooth solutions, the LISL scheme yields
a numerical order of convergence of one, whereas the MPCSL scheme
yields order two. The CPU time is now proportional to $\frac1{\Dx^3}$,
reflecting that we use $\frac{4\pi}{\Dx}$ grid points to discretize the control.
\begin{table}[ht]
\newsubtab{
3.93e-2&3.01e-2&&2.00&8.40e-4&&4.74\\
1.96e-2&1.61e-2&0.91&23.22&2.12e-4&1.98&53.06\\
9.82e-3&8.03e-3&1.00&268.64&5.30e-5&2.00&743.35\\
4.91e-3&3.94e-3&1.03&2161.63&1.33e-5&2.00&5995.62\\
2.45e-3&2.03e-3&0.96&17366.16&3.32e-6&2.00&48150.42}
{}{Tab:Ex3a}
\\\newsubtab{
3.93e-2&2.18e-2&&4.31&5.14e-4&&9.51\\
1.96e-2&1.07e-2&1.03&49.49&1.29e-4&2.00&109.50\\
9.82e-3&5.45e-3&0.97&571.99&3.21e-5&2.00&1515.24\\
4.91e-3&2.55e-3&1.10&4608.88&8.03e-6&2.00&12241.84\\
2.45e-3&1.34e-3&0.92&36995.81&2.01e-6&2.00&98378.95}
{}{Tab:Ex3b}
\caption{Results for optimal control problems at $t=0.5$, grid adapted to monotonicity\label{Tab:Ex3aundb}}
\end{table}
\subsection{Convergence test for a super-replication problem}\label{subsec:Bruder}
We consider a test problem from \cite{BBMZ} which was used to test
convergence rates for numerical approximations of a super-replication
problem from finance. The corresponding PDE is
\begin{equation}\label{eq:Superrepeq}
\inf_{\alpha_1^2+\alpha_2^2=1}
\left\{
\alpha_1^2u_t(t,x)-\frac12\tr\left(\sigma^\alpha(t,x)\sigma^{\alpha\,\top}(t,x) D^2u(t,x)\right)
\right\}
=f(t,x),\quad 0\leq x_1,x_2\leq 3
\end{equation}
with $\sigma^\alpha(t,x)=\begin{pmatrix}
\alpha_1x_1\sqrt{x_2}\\
\alpha_2\eta(x_2)
\end{pmatrix}$ and $\eta(x)=x(3-x)$. We take
$u(t,x)=1+t^2-e^{-x_1^2-x_2^2}$ as exact solution as in \cite{BBMZ},
and then $f$ is forced to be
\begin{align*}
&f(t,x)=
\frac12
\left(
u_t-\frac12x_1^2x_2u_{x_1x_1}-\frac12x_2^2(3-x_2)^2u_{x_2x_2}
\right.\\
&
\qquad
\left.
-\sqrt{\left(-u_t+\frac12x_1^2x_2u_{x_1x_1}-\frac12x_2^2(3-x_2)^2u_{x_2x_2}\right)^2+
\left(x_1\sqrt{x_2}^3(3-x_2)u_{x_1x_2}\right)^2
}
\right).
\end{align*}
In \cite{BBMZ} $\eta(x)=x$, while we take $\eta(x)=x(3-x)$
to prevent the LISL scheme from overstepping the boundaries.
Note that changing $\eta$ does {\em not} change the solutions as long
as $\eta>0$ in the interior of the domain, see \cite{BBMZ}, and hence
the above equation is equivalent to the equation used in \cite{BBMZ}.
The initial values and Dirichlet boundary values at $x_1=0$ and $x_2=0$ are
taken from the exact solution.  As
in \cite{BBMZ}, at $x=3$ and $y=3$ homogeneous Neumann boundary
conditions are implemented.  To approximate the values of
$\alpha_1,\alpha_2$, the Howard algorithm is used (see \cite{BBMZ}),
which requires an implicit time discretization, so we choose
$\theta=1$. As stop criterion of the iterations we require that the change of the maximal component and the sum over all components of the residual in Howard's algorithm are both smaller than $0.01$. The minimization is done over $\alpha_{1,k}+i\alpha_{2,k}=e^{2\pi ik/2N_{\dx}}$, $k=1,\dots,N_{\dx}$, where $N_{\dx}=3/\dx$ is the number of space grid points in one dimension. The linear systems involved are solved by the standard MATLAB back slash operator, using internally UMFPACK \cite{DA}. The numbers of time steps are chosen as $\frac1{\Dx}$ for the LISL scheme and $\frac1{\Dx^2}$ for the MPCSL scheme, respectively.

The results at $t=1$ are given in Table \ref{Tab:Ex7bundbc}.  Again, the
numerical order of convergence is approximately one when the LISL
scheme is used and approximately two for the MPCSL scheme. The CPU
times are better than expected for both the LISL and MPCSL schemes:
They get multiplied roughly by $10$ when $\Dx$
is divided by 2, a property which can also be observed in
\cite{BBMZ}. The reason is that the Howard algorithm needs fewer
iterations when the time step becomes smaller.

\begin{table}[ht]\captionsetup[subfloat]{position=top}
\subfloat[LISL]{\errtimesubtab{
1.50e-1&2.01e-1&1.80e-1&3.12e-1&&&& 0.71\\
7.50e-2&9.49e-2&9.58e-2&1.73e-1&1.08&0.91&0.85& 5.52\\
3.75e-2&4.29e-2&4.49e-2&8.27e-2&1.15&1.09&1.07&59.32\\
1.87e-2&1.94e-2&2.08e-2&3.88e-2&1.15&1.11&1.09&803.26}{Convergence results for Bruder type example with mixed Dirichlet and Neumann bc's when Howard's algorithm and linear interpolation are applied}{Tab:Ex7b}}
\subfloat[MPCSL]{\errtimesubtab{
3.00e-1&8.21e-2&7.65e-2&1.24e-1&&&& 1.40\\
1.50e-1&1.83e-2&1.81e-2&3.31e-2&2.17&2.08&1.90&11.38\\
7.50e-2&5.03e-3&4.10e-3&8.13e-3&1.86&2.14&2.03&124.25}{Convergence results for Bruder type example with mixed Dirichlet and Neumann bc's when Howard's algorithm and cubic interpolation are applied}{Tab:Ex7bc}}
\caption{Results for the convergence test for the super-replication problem at $t=1$}\label{Tab:Ex7bundbc}
\end{table}
\begin{remark}
Equation \eqref{eq:Superrepeq} can not be written in a form
\eqref{E} satisfying the assumptions of this paper, so the results of
this paper do not apply to this problem. However, it seems possible to
extend the analysis of Section \ref{Sec:Anal} to cover this
problem using comparison results from \cite{BBMZ} along with
$L^\infty$-bounds on the numerical solution that follow
from the maximum principle. Because of the unusual structure of
\eqref{eq:Superrepeq}, this analysis is not standard and outside
the scope of this paper.
\end{remark}
\begin{remark}
If we compare naively these results to the results of  \cite{BBMZ}, we
find that the LISL and MPCSL schemes are about 10 and up to a 1000 times faster
than the method of \cite{BBMZ}. Of course, this comparison is not
fair, e.\,g.\ it could be that a less efficient linear solver is used in
\cite{BBMZ}.
\end{remark}

\appendix

\section{Monotonicity preserving cubic interpolation}
\label{Sec:MPCI}

To define this type of interpolation, we start by
a 1D function $\phi$. For each sub-interval
$[x_i,x_{i+1}]$, $i\in\Z$, we construct a
cubic Hermite interpolant
\[\I\phi(x)=c_0+c_1(x-x_i)+c_2(x-x_i)^2+c_3(x-x_i)^3\]
fulfilling
\begin{align*}
\I\phi(x_i)&=\phi_i,&(\I\phi)'(x_i)&=d_i,&
\I\phi(x_{i+1})&=\phi_{i+1},&(\I\phi)'(x_{i+1})&=d_{i+1},\!\!\!\!
\end{align*}
where $\phi_i=\phi(x_i)$ and $d_i$ is an estimate of the
derivative of $\phi$ at $x_i$. It follows that
\begin{align}\label{Eq:CoeffDef}
c_0&=\phi_i,&c_1&=d_i,&
c_2&=\frac{3\Delta_i-d_{i+1}-2d_i}\dx,
&c_3&=-\frac{2\Delta_i-d_{i+1}-d_i}{\dx^2},
\end{align}
where $\Delta_i=\frac{\phi_{i+1}-\phi_i}\dx$.
To get a fourth order accurate interpolant,
$\phi_i'$ must be at least third order accurate, and we take
the symmetric fourth order approximation
\begin{equation}\label{Eq:DerAppr}
d_i=\frac{\phi_{i-2}-8\phi_{i-1}+8\phi_{i+1}-\phi_{i+2}}{12\dx},\quad i\in\Z.
\end{equation}
The resulting interpolation is not monotonicity preserving. Necessary and
sufficient conditions for preserving monotonicity were found by Fritsch and
Carlson \cite{FC} (see also \cite{RW}): If
$\Delta_i=0$, then monotonicity follows if and only if
$d_i=d_{i+1}=0$, and if
\[\alpha_i=\frac{d_i}{\Delta_i}\qquad\text{and}\qquad
\beta_i=\frac{d_{i+1}}{\Delta_i},\]
then monotonicity for  $\Delta_i\neq0$ follows if and only if
$(\alpha_i,\beta_i)\in\M=\M_e\cup\M_b$ where
\begin{align*}
M_e&=\{(\alpha,\beta):~(\alpha-1)^2+(\alpha-1)(\beta-1)+(\beta-1)^2-3(\alpha+\beta-2)\leq0\},\\
M_b&=\{(\alpha,\beta):~0\leq\alpha\leq3,~0\leq\beta\leq3\}.
\end{align*}
Eisenstat, Jackson and Lewis \cite{EJL} give an
algorithm that modifies the derivative approximation $d_i$ such that the
above conditions are fulfilled, and for monotone data the resulting
interpolant is a $C^1$ fourth order approximation. We will only
consider $C^0$ interpolants, and in that case their algorithm simplifies to the following steps to compute $(\I\phi)(x)$ on the interval $[x_i,x_{i+1}]$:

\smallskip

\begin{enumerate}[Step 1]
\medskip
\item Compute the initial $d_i$
using \eqref{Eq:DerAppr}.
\medskip
\item Compute $\Delta_i$. If $\Delta_i\neq0$ compute $\alpha_i$ and
  $\beta_i$, else set $\alpha_i=\beta_i=1$.
\medskip
\item Set $\alpha_i:=\max\{\alpha_i,0\}$ and $\beta_i:=\max\{\beta_i,0\}$.
\medskip
\item If $(\alpha_i,\beta_i)\notin\M$, modify $(\alpha_i,\beta_i)$ as follows:
\begin{itemize}\itemindent-3ex
\item If $\alpha_i\geq3$ and $\beta_i\geq3$, set $\alpha_i=\beta_i=3$,
\item else if $\beta_i>3$ and $\alpha_i+\beta_i\geq4$, decrease
$\beta_i$ such that $(\alpha_i,\beta_i)\in\partial\M$,
\item else if $\beta_i>3$ and $\alpha_i+\beta_i<4$, increase $\alpha_i$ such that $(\alpha_i,\beta_i)\in\partial\M$ or $\alpha_i=4-\beta_i$, in the last case subsequently decrease $\beta_i$ until $(\alpha_i,\beta_i)\in\partial\M$,
\item else if $\alpha_i>3$ and $\alpha_i+\beta_i\geq4$, decrease
$\alpha_i$ such that $(\alpha_i,\beta_i)\in\partial\M$,
\item else if $\alpha_i>3$ and $\alpha_i+\beta_i<4$, increase
  $\beta_i$ such that $(\alpha_i,\beta_i)\in\partial\M$ or
  $\beta_i=4-\alpha_i$, in the last case subsequently decrease
  $\alpha_i$ until $(\alpha_i,\beta_i)\in\partial\M$.
\end{itemize}
\medskip
\item Finally, replace $d_i$ by $\alpha_i\Delta_i$ and $d_{i+1}$ by
  $\beta_i\Delta_i$ in \eqref{Eq:CoeffDef} and compute $\I\phi$ from
  the resulting formula which then equals \eqref{MPCI}.
\end{enumerate}

\section{The proof of Theorem \ref{LipCont}}
\label{Sec:Pf}
We will prove the result when $k_0=0$. The general case can be reduced
to this case in a standard way by considering $U/R_{k_0}$ and
$\tU/R_{k_0}$ instead of $U$ and $\tU$.
We use doubling of variables techniques similar to those
used to prove this type of results for equation \eqref{E}.
 We take
\begin{align*}
&m_0=|(U(0,\cdot)-\tU(0,\cdot))^+|_0,\quad m=\sup_{\alp}\Big[|(f^\alp-\tf^\alp)^+|_0+(|U|_0\wedge|
    \tU|_0)|c^\alp-\tc^\alp|_0\Big],\\
& M^2= 4\sup_{\alp}
    \big[|\sigma^\alp-\ts^\alp|_0^2+|b^\alp-\tb^\alp|_0^2\big],
\end{align*}
where $\phi^+$ denotes the positive part of $\phi$,
and define $W(t,x,y)=U(t,x)-\tU(t,y)$,
\begin{align*}
& \phi(t,x,y)=m_0+tm+
\frac{1}{2\eps}K_T t
    M^2+\frac12R_{k_1}(t)(L_0+tL)(\eps
     +\frac1{\eps}|x-y|^2)+\delta(|x|^2+|y|^2),\\
    &\psi(t,x,y)=W(t,x,y) -
    \phi(t,x,y)- \eta (1+t),\quad\tilde m  =
    \sup_{\substack{t\in\dt\,\N_0\\x,y\in
        \R^N}}\psi(t,x,y) = \psi(\tilde t,\tx,\ty),
\end{align*}
for $\eps,\delta,\eta>0$ and a maximum point $(\tilde t,\tx,\ty)$. A
maximum point exists because of the $\delta$-terms in $\phi$.
We will prove that for any sequence $\eta_l\ra 0$, there is another
sequence $\delta_l\ra0$ such that $\psi(\tilde t_l,\tilde
x_l, \tilde y_l)\leq o(1)$ as $l\rightarrow\infty$. This
implies Theorem \ref{LipCont} when $k_0=0$. To see
this, fix $t>0,x,y$ and note that for any $\eps>0$,
\begin{align*}
&U(t,x)-\tU(t,y)-m_0-tm-  \frac{1}{2\eps} K_T t
     M^2-\frac12R_{k_1}(t)(L_0+tL)(\eps
     +\frac1{\eps}|x-y|^2)\\
&\leq \psi(\tilde t_l,\tilde x_l,
\tilde y_l)+\delta_l(|x|^2+|y|^2)+ \eta_l (1+t)\leq o(1) \quad\text{as}\quad
l\ra\infty.
\end{align*}
In this inequality we send $l\ra\infty$ and choose $\eps=|x-y|\vee t^{1/2}M$
to find that
\[U(t,x)-\tU(t,y)\leq m_0 + tm+t^{1/2}K_TM+R_{k_1}(t)(L_0+tL)|x-y|,\]
and hence Theorem \ref{LipCont} follows since $t>0,x,y$ were arbitrary.
We will not be explicit about the form of the
$\delta$-terms below. Their role is only to guarantee that the maximum is
attained at a (finite) point $(\tilde t,\tx,\ty)$, and their
contribution will always be $o(1)$ as $\delta\ra0$ (see also Section 3
in \cite{BJ:Rate}).

It is enough to prove that for every $\eta>0$, $\psi(\tilde
t,\tx,\ty)\leq o(1)$ as $\delta\ra0$. We proceed by contradiction
assuming there is an $\eta>0$ such that $\lim_{\delta\ra0}\psi(\tilde
t,\tx,\ty)> 0$. By the definition of $\psi$ we now
have $W(\tilde t,\tx,\ty)>0$ and $\tilde t>0$ for all $\delta>0$ small enough.
The last statement is true since
\[\psi(0,\tx,\ty)
\leq  m_0+L_0|\tx-\ty|-m_0-\frac {L_0}2(\eps
+\frac1{\eps}|\tx-\ty|^2)-\eta< 0.\]
The rest of the proof will aim at getting a contradiction for the case
$\tilde t>0$. Even if we do not write it like that, what we show below is
that $\frac{\psi(\tilde t, \tx,\ty)-\psi(\tilde t-\dt, \tx,\ty)}{\dt} \leq
o(1)-\eta$ as $\delta\ra0$, and this is impossible since $(\tilde
t,\tx,\ty)$ is a maximum point of $\psi$.

We proceed by defining the operator $\Pi^\alp$,
\begin{multline*}
\Pi^\alp
[\phi(t,\cdot,\cdot)](r,x,y)= \sum_{i=1}^M\Big\{\phi(t,x+y_{k,i}^{\alp,+}(r,x),y
+\tilde y_{k,i}^{\alp,+}(r,y))\\
-2\phi(t,x,y)+
\phi(t,x+y_{k,i}^{\alp,-}(r, x),y+\tilde y_{k,i}^{\alp,-}(r,y))\Big\}.
\end{multline*}
By the definition of $L^\alp_k$ and $\tL^\alp_k$,
it follows that
\begin{align*}
&\Pi^\alp[W(t,\cdot,\cdot)](r,x,y)= 2k^2\Big\{
L^\alp_k[U(t,\cdot)](r,x) - \tL^\alp_k[\tU(t,\cdot)](r,y)\Big\}.
\end{align*}
We set $\lambda:=\frac {\dt}{k^2}$ and
subtract the inequalities defining $U$ and $\tU$ (see
\eqref{FD1} and \eqref{FD2}) to find that for $(t,x),(t,y)\in Q_T$
\begin{multline*}
W(t,x,y)\leq  W(t-\dt,x,y)+\sup_\alp\Big\{\frac\lambda2 \Pi^\alp [
  \mm{W}(t,\cdot,\cdot)](t^\theta,x,y)\\
+ \dt\, c^\alp(t^\theta,x) \mm{W}(t,x,y)\Big\}
+\dt\,L|x-y|+ \dt\, m,
\end{multline*}
where $\mm{W}(t,x,y)=(1-\theta){W}(t-\dt,x,y)+\theta{W}(t,x,y)$ and
$t^\theta=t-(1-\theta)\dt$. Note that this new ``scheme'' is still
monotone by the definition of $\Pi^\alp$ and the
CFL condition. Hence we may replace $W$ in the above inequality by
any bigger function
coinciding with $W$ at $(t,x,y)$. By the definition of $\tilde m$,
\[W\leq \phi+\eta(1+t)+\tilde m\quad\text{in}\quad \dt\,\N_0\times
\R^N\times \R^N,\]
and equality holds
at $(\tilde t,\tx,\ty)$. Therefore we find that
\begin{multline}\tag{$\star$}
\label{phi-eq}
\phi(\tilde t,\tx,\ty)+\eta(1+\tilde t)\leq  \phi(\tilde
  t-\dt,\tx,\ty)+\eta(1+\tilde t-\dt)\\
+\sup_\alp\frac\lambda2\Pi^\alp [\mm{\phi}(\tilde t,\cdot,\cdot)](\tilde
t^\theta,\tx,\ty)\nonumber
+ \dt\, L|\tx-\ty|+ \dt\, m.\nonumber
\end{multline}
Here we also used the fact that $\Pi^\alp [\eta(1+t)+\tilde m]=0$ and
$c^\alp\leq0$. Moreover we can Taylor expand to see that
\begin{align*}
  &\Pi^\alp
  [\phi(t,\cdot,\cdot)](r,x,y)=\sum_{i=1}^M\Big\{(Y_i^++Y_i^-)\cdot
  D_x\phi+(\tY_i^++\tY_i^-)\cdot D_y\phi\nonumber\\
  &\quad+\frac12\tr[D_{xx}^2\phi\cdot( Y_i^+Y_i^{+\,\top}+
  Y_i^-Y_i^{-\,\top})]+\frac12\tr[D_{yy}^2\phi\cdot
  (\tY_i^+\tY_i^{+\,\top}+ \tY_i^-\tY_i^{-\,\top})]\nonumber\\
  &\quad+\frac12\tr[D_{xy}^2\phi\cdot (Y_i^+\tilde Y_i^{+\,\top}+
  \tilde Y_i^+ Y_i^{+\,\top}+Y_i^-\tilde Y_i^{-\,\top}+ \tilde Y_i^-
  Y_i^{-\,\top})]\Big\},\nonumber
\end{align*}
where $Y_i^\pm=y_{k,i}^{\alp,\pm}(r, x)$ and $\tilde Y_i^\pm=\tilde
y_{k,i}^{\alp,\pm}(r, y)$.
Now we use \eqref{eq:Y3} along with the definition of
$\phi$, to see that
\begin{align*}
&\Pi^\alp [\phi(t,\cdot,\cdot)](r,x,y)\leq  \frac1{\eps}
R_{k_1}(t)(L_0+tL) \bigg\{2k^2(b^\alp(r,x)-\tb^\alp(r,y))(x-y)\\
&\quad +
    k^2\tr\Big[(\sigma^\alp(r,x)-\tilde\sigma^\alp(r,y))(\sigma^\alp(r,x)
      -\tilde\sigma^\alp(r,y))^\top\Big]\\
&\quad + k^4\tr\Big[(b^\alp(r,x)-\tb^\alp(r,y))(b^\alp(r,x)
      -\tb^\alp(r,y))^\top\Big]\bigg\} + o(1),
\end{align*}
as $\delta\ra0$. These considerations lead to the following simplification of
\eqref{phi-eq},
\begin{align*}
&\eta  + \frac{\phi(\tilde t,\tx,\ty) - \phi(\tilde  t-\dt,\tx,\ty)} \dt\\
&\leq \theta \frac1{\eps} R_{k_1}(\tilde t)(L_0+\tilde tL)(\frac12M^2+\frac14 k_1|\tx-\ty|^2) \\
&\quad+ (1-\theta)
\frac1{\eps} R_{k_1}(\tilde t-\dt)(L_0+(\tilde t-\dt)L)(\frac12M^2+
\frac14k_1|\tx-\ty|^2)\\
&\quad + L|\tx-\ty|+m+o(1) \\
&\leq \frac1{\eps}R_{k_1}(\tilde t)(L_0+\tilde tL)(\frac12M^2+
\frac14 k_1|\tx-\ty|^2) +L|\tx-\ty|+ m+o(1):=RHS,
\end{align*}
as $\delta\ra0$. Now we proceed to calculate
$\delta_{\dt}\phi(t,x,y)=\frac{\phi(t,x,y) - \phi(t-\dt,x,y)}\dt.$
To do that we note that
$\delta_{\dt}(uv)=(\delta_{\dt} u)v+u\delta_{\dt} v-\dt (\delta_{\dt}
u)(\delta_{\dt} v).$
Since $\delta_{\dt}R_{k_1}(t)= k_1R_{k_1}(t)$ we then see that
\[\delta_{\dt}[R_{k_1}(t)(L_0+tL)]=k_1R_{k_1}(t)(L_0+tL)+R_{k_1}(t)L-\dt
Lk_1R_{k_1}(t),\]
and hence
\begin{align*}
&\delta_{\dt}\phi(\tilde t,\tx,\ty)
=m+\frac12 K_T \frac1\eps M^2+\frac12 R_{k_1}(\tilde t)\big[k_1(L_0+(\tilde t-\dt)L)+
L\big](\eps+\frac1\eps|\tx-\ty|^2).
\end{align*}
All of this leads to
$\eta\leq  RHS - \delta_{\dt}\phi(\tilde t,\tx,\ty)\leq
o(1)$ as
 $\delta\ra0$ and $2\Dt L\leq L_0$.
The last inequality follows from the bound on $K_T$. We have our
contradiction and the proof is complete.
\bibliographystyle{amsplain}

\begin{thebibliography}{10}

\bibitem{BJ:Rate}
Guy Barles and Espen~Robstad Jakobsen, \emph{On the convergence rate of
  approximation schemes for {H}amilton-{J}acobi-{B}ellman equations}, M2AN
  Math. Model. Numer. Anal. \textbf{36} (2002), no.~1, 33--54. \MR{1916291
  (2003h:65142)}

\bibitem{BJ:Rate3}
\bysame, \emph{Error bounds for monotone approximation schemes for parabolic
  {H}amilton-{J}acobi-{B}ellman equations}, Math. Comp. \textbf{76} (2007),
  no.~260, 1861--1893 (electronic). \MR{2336272 (2008i:65161)}

\bibitem{BS:Conv}
Guy Barles and Panagiotis~E. Souganidis, \emph{Convergence of approximation
  schemes for fully nonlinear second order equations}, Asymptotic Anal.
  \textbf{4} (1991), no.~3, 271--283. \MR{1115933 (92d:35137)}

\bibitem{BBMZ}
Olivier Bokanowski, Benjamin Bruder, Stefania Maroso, and Hasnaa Zidani,
  \emph{Numerical approximation for a superreplication problem under gamma
  constraints}, SIAM J. Numer. Anal. \textbf{47} (2009), no.~3, 2289--2320.
  \MR{2519604 (2010j:91213)}

\bibitem{BOZ}
Joseph~Fr\'{e}d\'{e}ric Bonnans, \'{E}lisabeth Ottenwaelter, and Hasnaa Zidani,
  \emph{A fast algorithm for the two dimensional {HJB} equation of stochastic
  control}, M2AN Math. Model. Numer. Anal. \textbf{38} (2004), no.~4, 723--735.
  \MR{2087732 (2005e:93165)}

\bibitem{BZ}
Joseph~Fr\'{e}d\'{e}ric Bonnans and Hasnaa Zidani, \emph{Consistency of
  generalized finite difference schemes for the stochastic {HJB} equation},
  SIAM J. Numer. Anal. \textbf{41} (2003), no.~3, 1008--1021. \MR{2005192
  (2004i:49061)}

\bibitem{CF:Appr}
Fabio Camilli and Maurizio Falcone, \emph{An approximation scheme for the
 optimal control of diffusion processes}, RAIRO Mod\'el. Math. Anal. Num\'er.
  \textbf{29} (1995), no.~1, 97--122. \MR{1326802 (96a:49033)}

\bibitem{CJ}
Fabio Camilli and Espen~Robstad Jakobsen, \emph{A finite element like scheme
  for integro-partial differential {H}amilton-{J}acobi-{B}ellman equations},
  SIAM J. Numer. Anal. \textbf{47} (2009), no.~4, 2407--2431. \MR{2525605
  (2010j:65227)}

\bibitem{CD}
Italo Capuzzo~Dolcetta, \emph{On a discrete approximation of the
  {H}amilton-{J}acobi equation of dynamic programming}, Appl. Math. Optim.
  \textbf{10} (1983), no.~4, 367--377. \MR{713483 (84j:49024)}

\bibitem{CIL:UG}
Michael~G. Crandall, Hitoshi Ishii, and Pierre-Louis Lions, \emph{User's guide
  to viscosity solutions of second order partial differential equations}, Bull.
  Amer. Math. Soc. (N.S.) \textbf{27} (1992), no.~1, 1--67. \MR{1118699
  (92j:35050)}

\bibitem{CL:MCM}
Michael~G. Crandall and Pierre-Louis Lions, \emph{Convergent difference schemes
  for nonlinear parabolic equations and mean curvature motion}, Numer. Math.
  \textbf{75} (1996), no.~1, 17--41. \MR{1417861 (97j:65134)}

\bibitem{crank47apm}
John Crank and Phyllis Nicolson, \emph{A practical method for numerical
  evaluation of solutions of partial differential equations of the
  heat-conduction type}, Proc. Cambridge Philos. Soc. \textbf{43} (1947),
  50--67. \MR{0019410 (8,409b)}

\bibitem{DA}
Timothy~A. Davis, \emph{Algorithm 832: {UMFPACK} {V}4.3---an
  unsymmetric-pattern multifrontal method}, ACM Trans. Math. Software
  \textbf{30} (2004), no.~2, 196--199. \MR{2075981}

\bibitem{DK:ConstCoeff}
Hongjie Dong and Nicolai~V. Krylov, \emph{On the rate of convergence of
  finite-difference approximations for {B}ellman equations with constant
  coefficients}, Algebra i Analiz \textbf{17} (2005), no.~2, 108--132.
  \MR{2159586 (2006f:49050)}

\bibitem{EJL}
Stanley~C. Eisenstat, Kenneth~R. Jackson, and John~W. Lewis, \emph{The order of
  monotone piecewise cubic interpolation}, SIAM J. Numer. Anal. \textbf{22}
  (1985), no.~6, 1220--1237. \MR{811195 (87d:65014)}

\bibitem{F}
Maurizio Falcone, \emph{A numerical approach to the infinite horizon problem of
  deterministic control theory}, Appl. Math. Optim. \textbf{15} (1987), no.~1,
  1--13. \MR{866164 (88c:49025)}

\bibitem{FC}
Frederic~N. Fritsch and Ralph~E. Carlson, \emph{Monotone piecewise cubic
  interpolation}, SIAM J. Numer. Anal. \textbf{17} (1980), no.~2, 238--246.
  \MR{567271 (81g:65012)}

\bibitem{Kr:00}
Nicolai~V. Krylov, \emph{On the rate of convergence of finite-difference
  approximations for {B}ellman's equations with variable coefficients}, Probab.
  Theory Related Fields \textbf{117} (2000), no.~1, 1--16. \MR{1759507
  (2001j:65134)}

\bibitem{KD:Book}
Harold~J. Kushner and Paul Dupuis, \emph{Numerical methods for stochastic
  control problems in continuous time}, second ed., Applications of Mathematics
  (New York), vol.~24, Springer-Verlag, New York, 2001, Stochastic Modelling
  and Applied Probability. \MR{1800098 (2001g:93002)}

\bibitem{M}
Jos\'{e}-Luis Menaldi, \emph{Some estimates for finite difference
  approximations}, SIAM J. Control Optim. \textbf{27} (1989), no.~3, 579--607.
  \MR{993288 (90m:65137)}

\bibitem{MW}
Theodore~Samuel Motzkin and Wolfgang Wasow, \emph{On the approximation of
  linear elliptic differential equations by difference equations with positive
  coefficients}, J. Math. Physics \textbf{31} (1953), 253--259. \MR{0052895
  (14,693i)}

\bibitem{MZ}
R\'{e}mi Munos and Hasnaa Zidani, \emph{Consistency of a simple
  multidimensional scheme for {H}amilton-{J}acobi-{B}ellman equations}, C. R.
  Math. Acad. Sci. Paris \textbf{340} (2005), no.~7, 499--502. \MR{2135230
  (2005k:49088)}

\bibitem{Ob}
Adam~M. Oberman, \emph{Convergent difference schemes for degenerate elliptic
  and parabolic equations: {H}amilton-{J}acobi equations and free boundary
  problems}, SIAM J. Numer. Anal. \textbf{44} (2006), no.~2, 879--895
  (electronic). \MR{2218974 (2007a:65173)}

\bibitem{Ob2}
\bysame, \emph{Wide stencil finite difference schemes for the elliptic
  {M}onge-{A}mp\`ere equation and functions of the eigenvalues of the
  {H}essian}, Discrete Contin. Dyn. Syst. Ser. B \textbf{10} (2008), no.~1,
  221--238. \MR{2399429 (2009f:35101)}

\bibitem{PFV}
David~M. Pooley, Peter~A. Forsyth, and Ken~R. Vetzal, \emph{Numerical
  convergence properties of option pricing {PDE}s with uncertain volatility},
  IMA J. Numer. Anal. \textbf{23} (2003), no.~2, 241--267. \MR{1974225
  (2004b:65132)}

\bibitem{RW}
Philip~J. Rasch and David~L. Williamson, \emph{On shape-preserving
  interpolation and semi-{L}agrangian transport}, SIAM J. Sci. Statist. Comput.
  \textbf{11} (1990), no.~4, 656--687. \MR{1054632 (91f:65020)}

\bibitem{YZ}
Jiongmin Yong and Xun~Yu Zhou, \emph{Stochastic controls}, Applications of
  Mathematics (New York), vol.~43, Springer-Verlag, New York, 1999, Hamiltonian
  systems and HJB equations. \MR{1696772 (2001d:93004)}

\end{thebibliography}
\providecommand{\bysame}{\leavevmode\hbox to3em{\hrulefill}\thinspace}
\providecommand{\MR}{\relax\ifhmode\unskip\space\fi MR }
\providecommand{\MRhref}[2]{%
  \href{http://www.ams.org/mathscinet-getitem?mr=#1}{#2}
}
\providecommand{\href}[2]{#2}

\end{document}